\documentclass[journal]{IEEEtran}
\usepackage[pdftex]{graphicx}
\usepackage{amsmath}
\usepackage{amssymb}
\usepackage{amsthm,subfigure}
\usepackage{color}
\usepackage{nccmath}
\usepackage{bbm}
\usepackage{hyperref}

\newtheorem{corollary}{Corollary}

\newtheorem{theorem}{Theorem}
\newtheorem{assumption}{Assumption}

\makeatletter
\newcommand{\pushright}[1]{\ifmeasuring@#1\else\omit\hfill$\displaystyle#1$\fi\ignorespaces}
\newcommand{\pushleft}[1]{\ifmeasuring@#1\else\omit$\displaystyle#1$\hfill\fi\ignorespaces}
\makeatother

\begin{document}

\title{Switching Control of Linear Time-Varying Networked Control Systems with Sparse Observer-Controller Networks}

\author{Mohammad~Razeghi-Jahromi$^{\dagger}$, Saeed Manaffam$^{\dagger\dagger}$, Alireza~Seyedi$^{\ddagger}$,~\IEEEmembership{Senior Member,~IEEE}, and Azadeh Vosoughi$^{\dagger\dagger}$,~\IEEEmembership{Member,~IEEE}
\thanks{$\dagger$M. Razeghi-Jahromi is with the ABB Corporate Research United States (USCRC),
Raleigh, NC 27606 USA (e-mail: mohammad.razeghi-jahromi@us.abb.com).}
\thanks{$\dagger\dagger$The authors are with the EECS Department, University of Central Florida, Orlando,
FL 32816 USA (e-mail: saeedmanaffam@knights.ucf.edu and azadeh@ucf.edu).}
\thanks{$\ddagger$A. Seyedi, deceased, was with the EECS Department, University of Central Florida, Orlando,
FL 32816 USA.}}

\maketitle
\begin{abstract}
In this paper we provide a set of stability conditions for linear time-varying networked control systems with arbitrary topologies using a piecewise quadratic switching stabilization approach with multiple quadratic Lyapunov functions. We use this set of stability conditions to provide a novel iterative low-complexity algorithm that must be updated and optimized in discrete time for the design of a sparse observer-controller network, for a given plant network with an arbitrary topology. We employ distributed observers by utilizing the output of other subsystems to improve the stability of each observer. To avoid unbounded growth of controller and observer gains, we impose bounds on the norms of the gains.
\end{abstract}



\section{Introduction}\label{sec:intro}

Networked control systems (NCS) have been the subject of much interest due to the fact that they have a wide range of applications, including electric power networks, transportation networks, factory automation, tele-operations and sensor and actuator networks, and the fact that they pose interesting unsolved problems in control theory. A centralized architecture is traditionally employed to control spatially distributed systems in which the components were connected via dedicated hard-wired links carrying the information from the sensors to a central location, where control signals were computed and sent to the actuators. However, the centralized architecture is not scalable. Moreover, it does not meet many new requirements such as modularity, resiliency, integrated diagnostics, and efficient maintenance. Distributed or decentralized networked control systems meet these requirements through reduction in the required communications and distribution of the computational power across the system. Consequently, these approaches can be scalable. More recently, the distributed and decentralized architectures have been made feasible due to relative maturity in communication and computing technologies, enabling their convergence with control.



\begin{figure}[!t]
\centering
\includegraphics[width =2.8in]{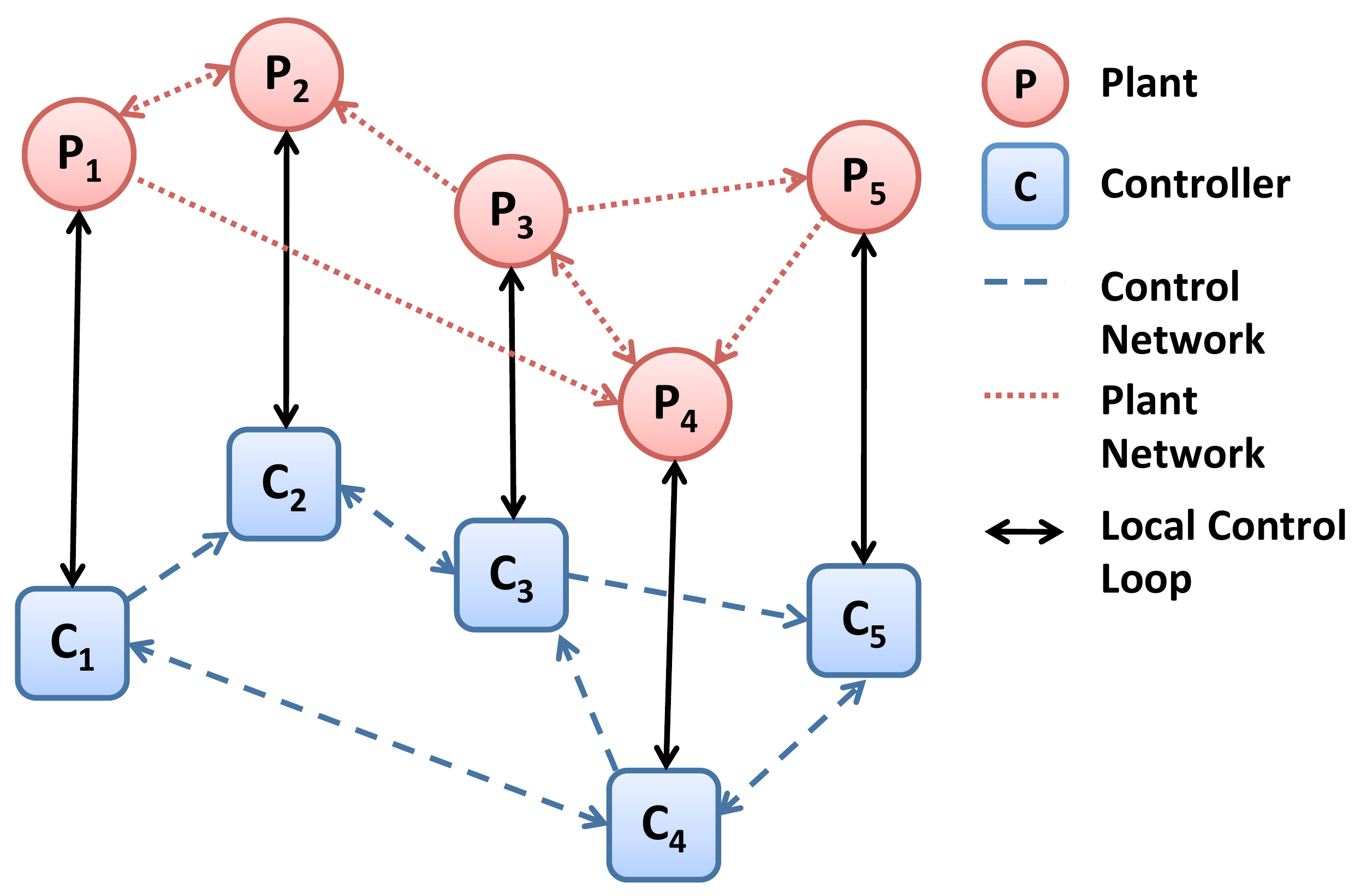}
\vspace{-0.1in}
\caption {A Networked Control System (NCS)}
\vspace{-0.2in}
\label{fig:NCS1}
\end{figure}

In general, a NCS consists of a number of subsystems, each comprised of a plant and a controller, coupled together in some network topology. The interaction of plants with one another forms the {\em plant network}. Measurements and control signals are communicated using the {\em control network}, a.k.a. information, communications, or feedback network (Fig. \ref{fig:NCS1}). This generalization covers the full range of architectures from decentralized, when the control network has no links, to centralized, when the control network is complete (i.e. all information is available to all controllers). For each architecture both the dynamics of each subsystem and the topology of plant network, play important roles in the stability of the entire interconnected system.



A key aspect in designing a particular NCS is the amount of information exchange. Typically, the all-to-all information exchange required for a centralized architecture is not feasible due to cost and complexity of the required communication. On the other end of the spectrum, ideally, one would have a decentralized controller \cite{Siljak1991}-\nocite{Siljak2005}\cite{Mazo2010}. It can be shown that even if each subsystem is asymptotically stable in isolation, the entire interconnected system may be unstable. Therefore, a decentralized architecture is usually inadequate to satisfy the performance requirements or even to stabilize the plant network. Thus, often the best solution is a distributed architecture, where some information is exchanged with neighbors \cite{Zhang2001}-\nocite{Baillieul2007,Hespanha2007,Yue2008}\cite{Wang2011}.

The networked control literature can be classified into two main groups. The first group as illustrated in \cite{montestruque2004}-\nocite{Wang2011,Teel2004,Walsh2002,seiler2005h,yu2005,wu2007,witrant2007,Gao08}\cite{Zhang13} studies different factors, including bandwidth, packet dropout and disorder, data quantization, time-varying sampling intervals, and time delays, all of which are imposed by imperfections and limitations of communication channels and can degrade the system performance or even destabilize the system. The second group, which this paper also falls into, considers the topological network effects, and investigates how the topology of the plant network affects the entire interconnected system stability and other performance metrics.


The problem of imposing {\em a priori} constraints on the controller has been arisen in previous articles on decentralized and distributed control of a general linear time-invariant (LTI) system. These constraints often are called the information constraints, and specify what information is available to which controller and manifests itself as sparsity or delay constraints.
In \cite{Lall2006}-\nocite{Rotkowitz2009,Rotkowitz2010}\cite{Swigart2009}, the authors have shown the convexification problem of finding optimal controller in order to minimize a norm of the closed loop map under a structural condition, namely {\em quadratic invariance}, which is an algebraic condition relating the plant and the constraint set. The works in \cite{Shah2011,Parrilo2011} have shown similar results, conditioned on the plant network being a partially ordered set (poset) which is based on poset information structure (acyclic information flow) among subsystems. This constraint is closely related to quadratic invariance, however, it can lead to more computationally efficient solutions and is applicable to a more general class of problems.


While these results are elegant and important, they impose restrictions on the topology of the plant network of LTI systems. The key question in the design of control network for NCSs with arbitrary plant network topologies is one of {\em topological information requirements} and can be framed as the question: {\em Which subsystems should be given the input and output information of a particular subsystem, in order for the local controllers to be able to satisfy a global control objective?} This is a critical question in the design of massively distributed control systems, such as the Smart Grid \cite{Massoud2011}-\nocite{Schuppen2011,Lafortune2003}\cite{Lin2007}.



In this paper, we extend our LTI NCS results \cite{Razeghi2015} by considering a linear time-varying (LTV) NCS with an arbitrary topology and provide a methodology to design an iterative sparse observer and controller network which updates in discrete time. As in our previous work \cite{Razeghi2015} here we also assume that the communication links do not have any bandwidth limitation, data loss or induced network delays.

We first use multiple quadratic Lyapunov functions to develop a set of stability conditions that guarantee global asymptotic stability using the piecewise quadratic switching stabilization method \cite{liberzon1999basic}-\nocite{liberzon2003switching,michel1999recent}\cite{lin2009stability}. We then use these stability conditions to design an iterative sparse observer-controller network for a given plant network with an arbitrary topology. We take a broader look at the topological information requirements by taking into account the distributed state estimation problem, which is generally neglected in existing work.


\section{Notation and Problem Definition}\label{sec:notation}

\subsubsection{Notation}
The set of real (column) $n$-vectors is denoted by $\mathbb{R}^{n}$ and the set of real $m\times n$ matrices is denoted by $\mathbb{R}^{m\times n}$. We use $\mathbb{R}_{+}$ and $\mathbb{Z}_{+}$ to denote the sets of non-negative real and non-negative integers respectively. Matrices and vectors are denoted by capital and lower-case bold letters, respectively. Generalized matrix inequality, $\prec$, is defined by the positive definite cone between symmetric matrices. The Euclidean ($l_2$) vector norm and the induced $l_2$ matrix norm are represented by $\| \cdot \|$ and the Frobenius matrix norm is denoted by $\|\cdot\|_F$. By $\lambda_{\mbox{\footnotesize min}}(\cdot)$, $\lambda_{\mbox{\footnotesize max}}(\cdot)$ and $\sigma_{\mbox{\footnotesize max}}(\cdot)$ we denote the smallest and largest eigenvalue and the largest singular value of the argument, respectively. The Schur (Hadamard) product is represented by $\circ$ and the $m\times n$ unit matrix consisting of all ones is denoted by $\mathbf{1}^{m\times n}$. We let $\mathcal{N}=\{1,\ldots,N\}$ and $\mathcal{N}_i = \mathcal{N}-\{i\}$. The indicator function of $x$ is represented by $1_{x}$ and column-stacking operator is denoted by $\mbox{vec}(\cdot)$.\

In the following subsection, we address the problem statement by employing the same methodology similar to that of \cite{Razeghi2015}.

\subsubsection{Problem Definition}
Consider a network of $N$ coupled LTV subsystems, each consisting of a plant and a controller. The state of the $i$th plant, $\mathbf{x}_i(t)\in \mathbb{R}^{n_i}$, is governed by
\begin{eqnarray}
\label{eq: 1}
\dot{\mathbf{x}}_i(t) &= & \textstyle{\mathbf{A}_i(t)\mathbf{x}_i(t)+\mathbf{B}_i(t)\mathbf{u}_i(t)+\sum_{j\in \mathcal{N}_i}\mathbf{H}_{ij}(t)\mathbf{x}_j(t)}\nonumber\\
\mathbf{y}_i(t) &= & \mathbf{C}_i(t)\mathbf{x}_i(t),
\end{eqnarray}
where $\mathbf{u}_i(t)\in \mathbb{R}^{m_i}$ and $\mathbf{y}_i(t)\in \mathbb{R}^{r_i}$ are input and output of the $i$th subsystem, and $\mathbf{A}_i(t)$, $\mathbf{B}_i(t)$, $\mathbf{C}_i(t)$ and $\mathbf{H}_{ij}(t)$ are known matrices. We assume that subsystem (\ref{eq: 1}) is both completely controllable and completely observable for all $i$. We consider an arbitrary directed network without self-loops. That is, $\mathbf{H}_{ii}(t)\equiv\mathbf{0}$, and $\mathbf{H}_{ij}(t)$ and $\mathbf{H}_{ji}(t)$ are not necessarily equal. We look for a distributed stabilizing observer-based controller of the form
\begin{eqnarray}
\label{eq: 2}
\dot{\hat{\mathbf{x}}}_i(t) &= & \textstyle{\mathbf{A}_i(t)\hat{\mathbf{x}}_i(t)+\mathbf{B}_i(t)\mathbf{u}_i(t)+\sum_{j\in \mathcal{N}_i}\mathbf{H}_{ij}(t)\hat{\mathbf{x}}_j(t)}\nonumber\\
&+&\textstyle{\mathbf{M}_i(t)(\mathbf{C}_i(t)\hat{\mathbf{x}}_i(t)-\mathbf{y}_i(t))}\nonumber\\
&+&\textstyle{\sum_{j\in \mathcal{N}_i}\mathbf{O}_{ij}(t)(\mathbf{C}_j(t)\hat{\mathbf{x}}_j(t)-\mathbf{y}_j(t))},\nonumber\\
\mathbf{u}_i(t) &= & \textstyle{\mathbf{K}_i(t)\hat{\mathbf{x}}_i(t)+\sum_{j\in \mathcal{N}_i}\mathbf{L}_{ij}(t)\hat{\mathbf{x}}_j(t)},
\end{eqnarray}
where $\hat{\mathbf{x}}_i(t)$ is the estimate of $\mathbf{x}_i(t)$, $\mathbf{K}_i(t)$ and $\mathbf{L}_{ij}(t)$ are local and coupling controller gains, and $\mathbf{M}_i(t)$ and $\mathbf{O}_{ij}(t)$ are local and coupling observer gains, respectively. Note that to estimate $\mathbf{x}_i(t)$, we not only use output of subsystem $i$, but also outputs of (potentially) all other subsystems. This is dual to the concept of distributed control. Our objective is to find distributed observer-based control law (\ref{eq: 2}), using feedback from (potentially) all other subsystems to stabilize the plant network with a sparse control network. That is, we aim to find $\mathbf{K}_i(t),\mathbf{M}_i(t),\mathbf{L}_{ij}(t)$ and $\mathbf{O}_{ij}(t)$, such that the overall network is globally asymptotically stable and that the number of links in the control network (number of non-zero coupling gains $\mathbf{L}_{ij}(t)$ and $\mathbf{O}_{ij}(t)$) is minimized. We also impose constraints
\begin{subequations}\label{eq: 5}
\begin{align}
\| \mathbf{K}_i(t)\| \leq &~ \kappa_i,\label{eq: 5a}\\
\| \mathbf{M}_i(t)\| \leq &~\mu_i,\label{eq: 5b}\\
\| \mathbf{L}_{ij}(t)\|\leq &~\iota_{ij},\label{eq: 5c}\\
\| \mathbf{O}_{ij}(t)\|\leq &~\omega_{ij},\label{eq: 5d}
\end{align}
\end{subequations}
to avoid undesirably large gains.

Defining $\mathbf{x}(t)=\mbox{vec}(\mathbf{x}_i(t))$, $\mathbf{u}(t)=\mbox{vec}(\mathbf{u}_i(t))$, $\mathbf{y}(t)=\mbox{vec}(\mathbf{y}_i(t))$, (\ref{eq: 1}) reduces to
\begin{align}
\label{eq: 1-1}
\dot{\mathbf{x}}(t) = &~ \mathbf{A}(t)\mathbf{x}(t)+\mathbf{B}(t)\mathbf{u}(t)+\mathbf{H}(t)\mathbf{x}(t),\nonumber\\
\mathbf{y}(t) = &~ \mathbf{C}(t)\mathbf{x}(t),
\end{align}
where $\mathbf{A}(t)=\mbox{diag}(\mathbf{A}_i(t))$, $\mathbf{B}(t)=\mbox{diag}(\mathbf{B}_i(t))$, $\mathbf{C}(t)=\mbox{diag}(\mathbf{C}_i(t))$ and $\mathbf{H}(t)=\left[\mathbf{H}_{ij}(t)\right]$.

Moreover, (\ref{eq: 2}) yields
\begin{eqnarray}
\label{eq: 1-3}
\dot{\hat{\mathbf{x}}}(t) &= & \mathbf{A}(t)\hat{\mathbf{x}}(t)+\mathbf{B}(t)\mathbf{u}(t)+\mathbf{H}(t)\hat{\mathbf{x}}(t)\nonumber\\
&+&\mathbf{M}(t)(\mathbf{C}(t)\hat{\mathbf{x}}(t)-\mathbf{y}(t))+\mathbf{O}(t)(\mathbf{C}(t)\hat{\mathbf{x}}(t)-\mathbf{y}(t)),\nonumber\\
\mathbf{u}(t) &= & \mathbf{K}(t)\hat{\mathbf{x}}(t)+\mathbf{L}(t)\hat{\mathbf{x}}(t),
\end{eqnarray}
where $\mathbf{K}(t)=\mbox{diag}(\mathbf{K}_i(t))$, $\mathbf{M}(t)=\mbox{diag}(\mathbf{M}_i(t))$, $\mathbf{L}(t)=\left[\mathbf{L}_{ij}(t)\right]$, with $\mathbf{L}_{ii}(t)\equiv\mathbf{0}$ and $\mathbf{O}(t)=\left[\mathbf{O}_{ij}(t)\right]$, with $\mathbf{O}_{ii}(t)\equiv\mathbf{0}$.\

Defining error $\mathbf{e}(t)\triangleq \hat{\mathbf{x}}(t)-\mathbf{x}(t)$ reduces (\ref{eq: 1-1}) and (\ref{eq: 1-3}) to
\begin{eqnarray}
\label{eq: 3}
\dot{\mathbf{x}}(t)&=&\left[\mathbf{A}(t)+\mathbf{H}(t)+\mathbf{B}(t)(\mathbf{K}(t)+\mathbf{L}(t))\right]\mathbf{x}(t)\nonumber\\
&+&\mathbf{B}(t)(\mathbf{K}(t)+\mathbf{L}(t))\mathbf{e}(t),\\
\label{eq: 4}
\dot{\mathbf{e}}(t)&=&\left[\mathbf{A}(t)+\mathbf{H}(t)+(\mathbf{M}(t)+\mathbf{O}(t))\mathbf{C}(t)\right]\mathbf{e}(t).
\end{eqnarray}
This is an LTV networked linear cascade dynamical system with the equilibrium point $(\mathbf{x},\mathbf{e})\equiv(\mathbf{0},\mathbf{0})$.


\begin{assumption} \label{assu:bounded} Matrices $\mathbf{A}(t), \mathbf{B}(t), \mathbf{C}(t)$ and $\mathbf{H}(t)$ are continuously differentiable and have bounded derivatives:
\begin{align}
\label{eq: A1}
\| \dot{\mathbf{A}}(t)\| \leq a,~\| \dot{\mathbf{B}}(t)\| \leq b,~\| \dot{\mathbf{C}}(t)\|\leq c,~\| \dot{\mathbf{H}}(t)\|\leq h.
\end{align}
\end{assumption}

Under Assumption \ref{assu:bounded}, for all $t, t'\in \mathbb{R}_{+}$ we have \cite{haddad2011}
\begin{align}
\label{eq: A2}
&\| \mathbf{A}(t)-\mathbf{A}(t')\| \leq a|t-t'|,~\| \mathbf{B}(t)-\mathbf{B}(t')\| \leq b|t-t'|,\nonumber\\
&\| \mathbf{C}(t)-\mathbf{C}(t')\|\leq c|t-t'|,~\| \mathbf{H}(t)-\mathbf{H}(t')\|\leq h|t-t'|.
\end{align}

\section{Piecewise Quadratic Switching Stabilization}
\label{sec:stability}

As the nature of problem is time-varying, and analytical design approaches are intractable (due to the generality of the network topology), one needs to take a numerical design approach. Thus, the design must be updated and optimized at discrete times. Consider the discrete series $0=t_0<t_1<\cdots<t_k<\cdots$. If a sample-and-hold control approach over intervals $[t_k,t_{k+1})$ is used, system (\ref{eq: 3}) and (\ref{eq: 4}) will become
\begin{eqnarray}
\label{eq: 3s}
\dot{\mathbf{x}}(t)&=&\left[\mathbf{A}(t)+\mathbf{H}(t)+\mathbf{B}(t)(\mathbf{K}_k+\mathbf{L}_k)\right]\mathbf{x}(t)\nonumber\\
&+&\mathbf{B}(t)(\mathbf{K}_k+\mathbf{L}_k)\mathbf{e}(t),\\
\label{eq: 4s}
\dot{\mathbf{e}}(t)&=&\left[\mathbf{A}(t)+\mathbf{H}(t)+(\mathbf{M}_k+\mathbf{O}_k)\mathbf{C}(t)\right]\mathbf{e}(t),
\end{eqnarray}
for $t\in[t_k,t_{k+1})$.



\begin{assumption}
For every $k\in \mathbb{Z}_{+}$, there exists constant $T_{\min}>0$ such that
\begin{align}
\label{eq: 3ss}
T_k\triangleq t_{k+1}-t_k\geq T_{\min}.
\end{align}
\end{assumption}

Later, in Corollary \ref{cor:Tmin} we will find a lower bound for $T_{\min}$.

\begin{theorem}
The equilibrium point, $\mathbf{x}\equiv\mathbf{0}$, of the system
\begin{align}
\label{eq: 3s-e}
\dot{\mathbf{x}}(t)=\left[\mathbf{A}(t)+\mathbf{H}(t)+\mathbf{B}(t)(\mathbf{K}_k+\mathbf{L}_k)\right]\mathbf{x}(t),
\end{align}
is globally asymptotically stable, if for every $k\in \mathbb{Z}_{+}$ and $t\in [t_k,t_{k+1})$, there exist $\mathbf{K}_k$, $\mathbf{L}_k$, $\mathbf{P}_k$, $\beta_i>0$, and $\epsilon_i>0$ such that
\begin{subequations}\label{eq: 7-e}
\begin{eqnarray}
\left[\mathbf{A}(t)+\mathbf{H}(t)+\mathbf{B}(t)(\mathbf{K}_k+\mathbf{L}_k)\right]^{T}\mathbf{P}_k+\quad\quad\,\,\nonumber \\
\mathbf{P}_k\left[\mathbf{A}(t)+\mathbf{H}(t)+\mathbf{B}(t)(\mathbf{K}_k+\mathbf{L}_k)\right]
+2\boldsymbol\beta\circ \mathbf{P}_k\preceq \mathbf{0}~~\label{eq: 7-ea}\\
\mathbf{P}_k-\boldsymbol\epsilon  \succeq \mathbf{0}~~\label{eq: 7-eb}\\
\mathbf{P}_{k-1}-\mathbf{P}_{k}\succeq \mathbf{0}~~\label{eq: 7-ec}
\end{eqnarray}
\end{subequations}
are satisfied for all $k\geq 1$ and \eqref{eq: 7-ea}, \eqref{eq: 7-eb} are satisfied for $k=0$, where $\mathbf{P}_k=\mbox{diag}(\mathbf{P}_{k,i})$, $\boldsymbol\beta=\mbox{diag}\left(\beta_i\mathbf{1}^{n_i\times n_i}\right)$, and $\boldsymbol\epsilon=\mbox{diag}(\epsilon_i \mathbf{I}_{n_i})$.
\end{theorem}

\begin{IEEEproof}
Consider multiple quadratic Lyapunov functions $V_k(\mathbf{x}(t))=\mathbf{x}^T(t)\mathbf{P}_k\mathbf{x}(t),~t\in [t_k,t_{k+1})$, and let $V(\mathbf{x}(t))=V_k(\mathbf{x}(t))$ when $t\in [t_k,t_{k+1})$. Since $\mathbf{P}_k\succ \mathbf{0}$, to show that $\mathbf{x}(t)\to 0$, it suffices to show that $V(\mathbf{x}(t))\to 0$.  Condition \eqref{eq: 7-ea} yields
\begin{align}
\label{eq: th1-e}
\dot{V}_k(\mathbf{x}(t))=&~\mathbf{x}^T(t)\left[\mathbf{A}(t)+\mathbf{H}(t)+\mathbf{B}(t)(\mathbf{K}_k+\mathbf{L}_k)\right]^{T}\mathbf{P}_k\mathbf{x}(t)&\nonumber\\
&~+\mathbf{x}^T(t)\mathbf{P}_k\left[\mathbf{A}(t)+\mathbf{H}(t)+\mathbf{B}(t)(\mathbf{K}_k+\mathbf{L}_k)\right]\mathbf{x}(t)\nonumber\\
\leq&~ -2\mathbf{x}^T(t)\boldsymbol\beta \circ \mathbf{P}_k\mathbf{x}(t)\nonumber\\
\leq&~ -2\beta\sum_{i\in \mathcal{N}}\mathbf{x}^T_i(t)\mathbf{P}_{k,i}\mathbf{x}_i(t)\nonumber\\
=&~-2\beta V_k(\mathbf{x}(t)),
\end{align}
for $t\in [t_k,t_{k+1})$, where $\beta=\min_i\beta_i$. This means that $\dot{V}(\mathbf{x}(t))=\dot{V}_k(\mathbf{x}(t))$ is negative definite over $t\in [t_k,t_{k+1})$. Moreover, since $\mathbf{x}(t)$ is continuous, condition \eqref{eq: 7-ec} implies that $V(\mathbf{x}(t))$ is decreasing over all $t$. Thus, since $V(\mathbf{x}(t))$ is positive definite, if we show that its samples, $V(\mathbf{x}(t_k))=V_k(\mathbf{x}(t_k))$, converge to zero we have $V(\mathbf{x}(t))\to 0$ as well. Since the sequence $V_k(\mathbf{x}(t_k))$ is monotonically decreasing and positive definite, it converges to some value $v\ge 0$. To show that $v=0$, we first note that due to the comparison lemma \cite{haddad2011}, \eqref{eq: th1-e} yields
\begin{align}
V_k(\mathbf{x}(t))\leq&~ V_k(\mathbf{x}(t_k))e^{-2\beta(t-t_k)},t\in [t_k,t_{k+1}).
\end{align}
Now, we can write
\begin{align}
\label{eq: th3}
0=&~v-v\nonumber\\
=&~\lim_{k\rightarrow \infty}V_{k+1}(\mathbf{x}(t_{k+1}))-\lim_{k\rightarrow \infty}V_{k}(\mathbf{x}(t_{k}))\nonumber\\
=&~\lim_{k\rightarrow \infty}\left[V_{k+1}(\mathbf{x}(t_{k+1}))-V_{k}(\mathbf{x}(t_{k+1}^{-}))\right]+\nonumber\\
&~\lim_{k\rightarrow \infty}\left[V_{k}(\mathbf{x}(t_{k+1}^{-}))-V_{k}(\mathbf{x}(t_{k}))\right]\nonumber\\
=&~\lim_{k\rightarrow \infty}\left[\mathbf{x}^T(t_{k+1})\left(\mathbf{P}_{k+1}-\mathbf{P}_k\right)\mathbf{x}(t_{k+1})\right]+\nonumber\\
&~\lim_{k\rightarrow \infty}\left[V_k(\mathbf{x}(t_k))\left(e^{-2\beta(t_{k+1}^{-}-t_k)}-1\right)\right]\nonumber\\
\leq &~\lim_{k\rightarrow \infty}\left[V_k(\mathbf{x}(t_k))\left(e^{-2\beta(t_{k+1}^{-}-t_k)}-1\right)\right]\nonumber\\
\leq&~ -\left(1-e^{-2\beta T_{\min}}\right)\lim_{k\rightarrow \infty}V_{k}(\mathbf{x}(t_k))\nonumber\\
\leq &~ -\left(1-e^{-2\beta T_{\min}}\right)\lim_{k\rightarrow \infty}\left[\lambda_{\mbox{\footnotesize min}}(\mathbf{P}_k)
\|\mathbf{x}(t_k)\|^2\right]\nonumber\\
\leq&~  -\epsilon\left(1-e^{-2\beta T_{\min}}\right)\lim_{k\rightarrow \infty}\|\mathbf{x}(t_k)\|^2\nonumber\\
\leq&~ 0,
\end{align}
where $\epsilon=\min_i\epsilon_i>0$. This implies that
\begin{align}
\label{eq: th4}
\epsilon\left(1-e^{-2\beta T_{\min}}\right)\lim_{k\rightarrow \infty}\|\mathbf{x}(t_k)\|^2=0.
\end{align}
which requires that $\lim_{k\rightarrow \infty} \mathbf{x}(t_k)\rightarrow 0$, and $\lim_{k\rightarrow \infty} V_k(\mathbf{x}(t_k))\rightarrow 0$. In other words, the system \eqref{eq: 3s-e} is globally asymptotically stable.
\end{IEEEproof}

\begin{theorem}\label{theononconvex}
The equilibrium point of the system in (\ref{eq: 3}) and (\ref{eq: 4}), $(\mathbf{x},\mathbf{e})\equiv(\mathbf{0},\mathbf{0})$, is globally asymptotically stable, if for every $k\in \mathbb{Z}_{+}$ and $t\in [t_k,t_{k+1})$, there exist $\mathbf{K}_k$, $\mathbf{L}_k$, $\mathbf{P}_k$, $\mathbf{M}_k$, $\mathbf{O}_k$, $\hat{\mathbf{P}}_k$, $\beta_i>0$, and $\epsilon_i>0$ such that
\begin{subequations}\label{eq: 7}
\begin{eqnarray}
\left[\mathbf{A}(t)+\mathbf{H}(t)+\mathbf{B}(t)(\mathbf{K}_k+\mathbf{L}_k)\right]^{T}\mathbf{P}_k+\quad\quad\,\,\nonumber \\
\mathbf{P}_k\left[\mathbf{A}(t)+\mathbf{H}(t)+\mathbf{B}(t)(\mathbf{K}_k+\mathbf{L}_k)\right]
+2\boldsymbol\beta\circ \mathbf{P}_k\preceq \mathbf{0}~~\label{eq: 7a}\\
\left[\mathbf{A}(t)+\mathbf{H}(t)+(\mathbf{M}_k+\mathbf{O}_k)\mathbf{C}(t)\right]^{T}\hat{\mathbf{P}}_k+\quad\quad\,\,\nonumber\\
\hat{\mathbf{P}}_k\left[\mathbf{A}(t)+\mathbf{H}(t)+(\mathbf{M}_k+\mathbf{O}_k)\mathbf{C}(t)\right]
+2\boldsymbol\beta\circ \hat{\mathbf{P}}_k\preceq \mathbf{0}~~\label{eq: 7b}\\
\mathbf{P}_k-\boldsymbol\epsilon  \succeq \mathbf{0}  ~~\label{eq: 7c}\\
\hat{\mathbf{P}}_k-\boldsymbol\epsilon  \succeq  \mathbf{0}  ~~\label{eq: 7d}\\
\mathbf{P}_{k-1}-\mathbf{P}_k\succeq \mathbf{0} ~~\label{eq: 7e}\\
\hat{\mathbf{P}}_{k-1}-\hat{\mathbf{P}}_k\succeq \mathbf{0}~~\label{eq: 7f}
\end{eqnarray}
\end{subequations}
are satisfied for all $k\ge 1$ and \eqref{eq: 7a} to \eqref{eq: 7d} are satisfied for $k=0$, where $\mathbf{P}_k=\mbox{diag}(\mathbf{P}_{k,i})$, $\hat{\mathbf{P}}_k=\mbox{diag}(\hat{\mathbf{P}}_{k,i})$, $\boldsymbol\beta=\mbox{diag}\left(\beta_i\mathbf{1}^{n_i\times n_i}\right)$, and $\boldsymbol\epsilon=\mbox{diag}(\epsilon_i \mathbf{I}_{n_i})$.
\end{theorem}

\begin{IEEEproof}
By Theorem 1, conditions \eqref{eq: 7a}, \eqref{eq: 7c} and \eqref{eq: 7e} imply that the unforced system \eqref{eq: 3s} is globally asymptotically stable and, consequently, \eqref{eq: 3s} is input-to-state stable. Similarly, \eqref{eq: 7b}, \eqref{eq: 7d} and \eqref{eq: 7f} guarantee global asymptotical stability of \eqref{eq: 4s}. Hence the equilibrium point of cascaded dynamical system \eqref{eq: 3s} and \eqref{eq: 4s}, $(\mathbf{x},\mathbf{e})\equiv(\mathbf{0},\mathbf{0})$, is globally asymptotically stable \cite{haddad2011}.
\end{IEEEproof}

We note that conditions \eqref{eq: 7a} and \eqref{eq: 7b} are the requirement that the energy in the system is decreasing in each interval while conditions \eqref{eq: 7c} and \eqref{eq: 7d} guarantee positive definiteness of multiple quadratic Lyapunov functions. Conditions \eqref{eq: 7e} and \eqref{eq: 7f} are necessary to guarantee multiple quadratic Lyapunov functions form a non-increasing sequence when entering the next interval.

Theorem \ref{theononconvex} provides a set of stability conditions based on which the controllers and observers can be designed. However, these conditions are not convex. This is important since a solution must be found iteratively at each $t_k$. The following theorem provides a set of convex, albeit more conservative, stability conditions. It also incorporates bounds \eqref{eq: 5}.

\begin{theorem}
\label{theo:Th3}
System (\ref{eq: 1}) with controller (\ref{eq: 2}) is globally asymptotically stable, and bounds (\ref{eq: 5}) are satisfied, if the following convex constraints have a feasible point, such that
\begin{subequations}\label{eq: 39}
\begin{eqnarray}
\mathbf{F}_k+\mathbf{F}^{T}_k+\gamma \mathbf{I}_{n}\preceq \mathbf{0}  \label{eq: 39a}\\
\hat{\mathbf{F}}_k+\hat{\mathbf{F}}^{T}_k+\gamma \mathbf{I}_{n}\preceq \mathbf{0}  \label{eq: 39b}\\
\boldsymbol\epsilon^{-1}\succeq\mathbf{Z}_k \succ \mathbf{0}  \label{eq: 39c}\\
\hat{\mathbf{P}}_k-\boldsymbol\epsilon \succeq \mathbf{0}  \label{eq: 39d}\\
\mathbf{Z}_{k}-\mathbf{Z}_{k-1} \succeq \mathbf{0}  \label{eq: 39e}\\
\hat{\mathbf{P}}_{k-1}-\hat{\mathbf{P}}_{k} \succeq \mathbf{0}  \label{eq: 39f}\\
\kappa_i\lambda_{\mbox{\footnotesize min}}(\mathbf{Z}_{k,i})-\sigma_{\mbox{\footnotesize max}}(\mathbf{W}_{k,i}) \geq  0  \label{eq: 39g}\\
\iota_{ij}\lambda_{\mbox{\footnotesize min}}(\mathbf{Z}_{k,j})-\sigma_{\mbox{\footnotesize max}}(\mathbf{Y}_{k,ij}) \geq  0  \label{eq: 39h}\\
\mu_i\lambda_{\mbox{\footnotesize min}}(\hat{\mathbf{P}}_{k,i})-\sigma_{\mbox{\footnotesize max}}(\hat{\mathbf{W}}_{k,i}) \geq 0  \label{eq: 39i}\\
\omega_{ij}\lambda_{\mbox{\footnotesize min}}(\hat{\mathbf{P}}_{k,i})-\sigma_{\mbox{\footnotesize max}}(\hat{\mathbf{Y}}_{k,ij}) \geq  0  \label{eq: 39j}
\end{eqnarray}
\end{subequations}
are satisfied for all $i,j\in \mathcal{N}$, $k\geq 1$, some $\gamma>0$ and \eqref{eq: 39a}-\eqref{eq: 39d} and \eqref{eq: 39g}-\eqref{eq: 39j} are satisfied for $k=0$, where
\begin{align}
\mathbf{F}_k=(\mathbf{A}_k+\mathbf{H}_k)\mathbf{Z}_k+\mathbf{B}_k(\mathbf{W}_k+\mathbf{Y}_k)
+\boldsymbol{\beta}\circ \mathbf{Z}_k,\nonumber\\
\hat{\mathbf{F}}_k=\hat{\mathbf{P}}_k(\mathbf{A}_k+\mathbf{H}_k)+(\hat{\mathbf{W}}_k
+\hat{\mathbf{Y}}_k)\mathbf{C}_k+\boldsymbol{\beta}\circ \hat{\mathbf{P}}_k,\nonumber
\end{align}
$\mathbf{Z}_k=\mbox{diag}(\mathbf{Z}_{k,i})$, $\hat{\mathbf{P}}_k=\mbox{diag}(\hat{\mathbf{P}}_{k,i})$, $\mathbf{W}_k=\mbox{diag}(\mathbf{W}_{k,i})$, $\hat{\mathbf{W}}_k=\mbox{diag}(\hat{\mathbf{W}}_{k,i})$, $\mathbf{Y}_k=[\mathbf{Y}_{k,ij}]$ and $\hat{\mathbf{Y}}_k=[\hat{\mathbf{Y}}_{k,ij}]$ with $\mathbf{Y}_{k,ii}=\hat{\mathbf{Y}}_{k,ii}=\mathbf{0}$, $\mathbf{A}_k=\mathbf{A}(t_k)$, $\mathbf{B}_k=\mathbf{B}(t_k)$, $\mathbf{C}_k=\mathbf{C}(t_k)$, and $\mathbf{H}_k=\mathbf{H}(t_k)$. \\
Furthermore, if $\mathbf{Z}^\star_{k,i},\hat{\mathbf{P}}^\star_{k,i},\mathbf{W}^\star_{k,i},\hat{\mathbf{W}}^\star_{k,i},\mathbf{Y}^\star_{k,ij}$, $\hat{\mathbf{Y}}^\star_{k,ij}$ is a solution of \eqref{eq: 39}, the controller and observer gains are
\begin{eqnarray}
\label{eq: 39.1}
\begin{array}{ll}
\mathbf{K}^\star_{k,i}=\mathbf{W}^\star_{k,i}\mathbf{Z}^{\star~-1}_{k,i},  & \mathbf{L}^\star_{k,ij}=\mathbf{Y}^{\star}_{k,ij}\mathbf{Z}^{\star~-1}_{k,j},\\
\mathbf{M}^\star_{k,i}=\hat{\mathbf{P}}^{\star~-1}_{k,i}\hat{\mathbf{W}}^\star_{k,i}, & \mathbf{O}^\star_{k,ij}=\hat{\mathbf{P}}^{\star~-1}_{k,i}\hat{\mathbf{Y}}^{\star}_{k,ij},
\end{array}
\end{eqnarray}
for all $t\in [t_k,t_{k+1})$ and the next switching time is $t_{k+1}=t_k+T_k$ where
\begin{align}
\label{eq: 39.14}
T_k=\frac{1}{2}\min\Bigg\{&\frac{-\lambda_{\mbox{\footnotesize max}}(\mathbf{F}_k+\mathbf{F}^T_k)}{(a+h)\|\mathbf{Z}_k\|+b\|\mathbf{W}_k+\mathbf{Y}_k\|},\nonumber\\
&\frac{-\lambda_{\mbox{\footnotesize max}}(\hat{\mathbf{F}}_k+\hat{\mathbf{F}}^T_k)}{(a+h)\|\hat{\mathbf{P}}_k\|+c\|\hat{\mathbf{W}}_k+\hat{\mathbf{Y}}_k\|}\Bigg\}.
\end{align}
\end{theorem}

\begin{IEEEproof}
We will show that the conditions of Theorem \ref{theononconvex}, namely \eqref{eq: 7a}-\eqref{eq: 7f}, are satisfied if \eqref{eq: 39a}-\eqref{eq: 39f} hold. Defining new variables $\mathbf{Z}_k\triangleq  \mathbf{P}^{-1}_k$, $\mathbf{W}_k\triangleq  \mathbf{K}_k\mathbf{P}^{-1}_k$, and $\mathbf{Y}_k\triangleq \mathbf{L}_k\mathbf{P}^{-1}_k$, reveals that \eqref{eq: 7a}, \eqref{eq: 7c}  and \eqref{eq: 7e} are equivalent to
\begin{subequations}\label{eq: 19}
\begin{eqnarray}
\left[(\mathbf{A}(t)+\mathbf{H}(t))\mathbf{Z}_k+\mathbf{B}(t)(\mathbf{W}_k+\mathbf{Y}_k)\right]^T+\quad\quad\,\,\,\,\,\nonumber\\
\left[(\mathbf{A}(t)+\mathbf{H}(t))\mathbf{Z}_k+\mathbf{B}(t)(\mathbf{W}_k+\mathbf{Y}_k)\right]
+2\boldsymbol \beta \circ \mathbf{Z}_k\preceq \mathbf{0}~~~\label{eq: 19a}\\
\boldsymbol\epsilon^{-1}\succeq\mathbf{Z}_k \succ \mathbf{0}~~~\label{eq: 19b}\\
\mathbf{Z}_{k}-\mathbf{Z}_{k-1}\succeq \mathbf{0}~~~\label{eq: 19c}
\end{eqnarray}
\end{subequations}
Clearly, \eqref{eq: 19b} and \eqref{eq: 19c} are \eqref{eq: 39c} and \eqref{eq: 39e}. To show that \eqref{eq: 19a} yields \eqref{eq: 39a}, we first note that
\begin{align}
(\mathbf{A}(t)+\mathbf{H}(t))\mathbf{Z}_k=&(\mathbf{A}_k+\mathbf{H}_k)\mathbf{Z}_k+(\bigtriangleup\mathbf{A}(t)+\bigtriangleup\mathbf{H}(t))\mathbf{Z}_k\nonumber\\
\preceq &(\mathbf{A}_k+\mathbf{H}_k)\mathbf{Z}_k+\nonumber\\
&\left\|(\bigtriangleup\mathbf{A}(t)+\bigtriangleup\mathbf{H}(t))\mathbf{Z}_k\right\|\mathbf{I}_n\nonumber\\
\preceq & (\mathbf{A}_k+\mathbf{H}_k)\mathbf{Z}_k+T_k(a+h)\|\mathbf{Z}_k\|\mathbf{I}_n,
\end{align}
where $\bigtriangleup\mathbf{A}(t)=\mathbf{A}(t)-\mathbf{A}_k$ and $\bigtriangleup\mathbf{H}(t)=\mathbf{H}(t)-\mathbf{H}_k$, and the last inequality is due to Assumption \ref{assu:bounded}. Similarly, we can show that
\begin{align}
\mathbf{B}(t)(\mathbf{W}_k+\mathbf{Y}_k)\preceq & \mathbf{B}_k(\mathbf{W}_k+\mathbf{Y}_k)+T_kb\|\mathbf{W}_k+\mathbf{Y}_k\|\mathbf{I}_n.
\end{align}
Thus, the left hand side of \eqref{eq: 19a} can be upper bounded by
\begin{eqnarray}
\mathbf{F}_k+\mathbf{F}^{T}_k+2T_k[(a+h)\|\mathbf{Z}_k\|+b\|\mathbf{W}_k+\mathbf{Y}_k\|]\mathbf{I}_n.\label{eq:ub}
\end{eqnarray}
In other words, if \eqref{eq:ub} is negative semidefinite, \eqref{eq: 7a} holds. Upper bound \eqref{eq:ub} is negative semidefinite, if there exist $\gamma >0$ such that
\begin{align}
\label{eq: LMI3}
\mathbf{F}_k+\mathbf{F}^T_k+\gamma \mathbf{I}_{n}\preceq \mathbf{0},
\end{align}
which is \eqref{eq: 39a} and
\begin{align}
\label{eq: T1}
T_k\le \frac{1}{2}\frac{-\lambda_{\mbox{\footnotesize max}}(\mathbf{F}_k+\mathbf{F}^T_k)}{(a+h)\|\mathbf{Z}_k\|+b\|\mathbf{W}_k+\mathbf{Y}_k\|},
\end{align}
which is guaranteed by \eqref{eq: 39.14}. A similar argument, omitted for brevity, shows that \eqref{eq: 7b}, \eqref{eq: 7d}  and \eqref{eq: 7f} are satisfied if \eqref{eq: 39b}, \eqref{eq: 39d}, \eqref{eq: 39f} and \eqref{eq: 39.14} hold.

For \eqref{eq: 5a} we note that we can upper bound the norm of $\mathbf{K}_i$ as
\begin{align}
\label{eq: 27}
\| \mathbf{K}_{k,i}\|  =&~  \| \mathbf{W}_{k,i}\mathbf{Z}^{-1}_{k,i}\| \nonumber\\
\leq &~ \| \mathbf{W}_{k,i}\|\| \mathbf{Z}^{-1}_{k,i}\|  \nonumber\\
=&~ \sigma_{\mbox{\footnotesize max}}(\mathbf{W}_{k,i})\lambda_{\mbox{\footnotesize max}}(\mathbf{Z}^{-1}_{k,i})\nonumber\\
=&~ \frac{\sigma_{\mbox{\footnotesize max}}(\mathbf{W}_{k,i})}{\lambda_{\mbox{\footnotesize min}}(\mathbf{Z}_{k,i})}.
\end{align}
Thus, forcing \eqref{eq: 5a} will be forced if
\begin{align}
\label{eq: 28}
\| \mathbf{K}_{k,i}\|\leq \frac{\sigma_{\mbox{\footnotesize max}}(\mathbf{W}_{k,i})}{\lambda_{\mbox{\footnotesize min}}(\mathbf{Z}_{k,i})}\leq \kappa_i,
\end{align}
or equivalently $\kappa_i\lambda_{\mbox{\footnotesize min}}(\mathbf{Z}_{k,i})-\sigma_{\mbox{\footnotesize max}}(\mathbf{W}_{k,i})\geq 0$, which is \eqref{eq: 39g}. Similarly, \eqref{eq: 5b}-\eqref{eq: 5d} are forced if \eqref{eq: 39h}-\eqref{eq: 39j} hold.

Finally, we note that the original variables can then be found from $\mathbf{P}_k=\mathbf{Z}^{-1}_k$, $\mathbf{K}_k=\mathbf{W}_k\mathbf{Z}^{-1}_k$ and $\mathbf{L}_k=\mathbf{Y}_k\mathbf{Z}^{-1}_k$.
\end{IEEEproof}

\begin{corollary}
\label{cor:Tmin}
If the conditions of Theorem 3 hold, a lower bound for $T_{\min}$ in Assumption 2 is
\begin{align}
\label{eq: T51}
T_{\min}=&\min_k T_k\nonumber\\
\geq&
\frac{1}{2}\min\Bigg\{\frac{\epsilon\gamma}{a+h+b\left(\kappa+\iota\right)},
~\frac{1}{\|\hat{\mathbf{P}}_0\|}\frac{\gamma}{a+h+c\left(\mu+\omega\right)}\Bigg\}
\end{align}
where $\epsilon=\min_i \epsilon_i>0$, $\gamma>0$ are the margins in inequalities \eqref{eq: 39a} and \eqref{eq: 39b} and $\kappa=\sum_{i\in \mathcal{N}} \sqrt{\min \{m_i,n_i\}}\kappa_i$, $\iota=\sum_{i,j\in \mathcal{N}}\sqrt{\min \{m_i,n_j\}}\iota_{ij}$, $\mu=\sum_{i\in \mathcal{N}} \sqrt{\min \{n_i,r_i\}}\mu_i$ and $\omega=\sum_{i,j\in \mathcal{N}}\sqrt{\min \{n_i,r_j\}}\omega_{ij}$.
\end{corollary}

\begin{IEEEproof}
Equation \eqref{eq: 39a} gives us $\lambda_{\mbox{\footnotesize max}}(\mathbf{F}_k+\mathbf{F}^T_k)\leq -\gamma$. From (\ref{eq: T1}), we have
\begin{align}
\label{eq: T11}
T_k=&~-\frac{1}{2}\frac{\lambda_{\mbox{\footnotesize max}}(\mathbf{F}_k+\mathbf{F}^T_k)}{(a+h)\|\mathbf{Z}_k\|+b\|\mathbf{W}_k+\mathbf{Y}_k\|}\nonumber\\
=&~-\frac{1}{2}\frac{\lambda_{\mbox{\footnotesize max}}(\mathbf{F}_k+\mathbf{F}^T_k)}{(a+h)\|\mathbf{Z}_k\|+b\|\mathbf{K}_k\mathbf{Z}_k+\mathbf{L}_k\mathbf{Z}_k\|}\nonumber\\
\geq &~\frac{1}{2}\frac{\epsilon\gamma}{a+h+b\left(\|\mathbf{K}_k\|_F+\|\mathbf{L}_k\|_F\right)}\nonumber\\
\geq &~\frac{1}{2}\frac{\epsilon\gamma}{a+h+b\left(\kappa+\iota\right)}.
\end{align}
The last four inequalities are satisfied because for any $\mathbf{D}\in \mathbb{R}^{m\times n}$, the Frobenius and Euclidian norms satisfy $\|\mathbf{D}\|\leq \|\mathbf{D}\|_F\leq \sqrt{\min \{m,n\}}\|\mathbf{D}\|$.

Similarly, we have
\begin{align}
\label{eq: T21}
T_k=&~-\frac{1}{2}\frac{\lambda_{\mbox{\footnotesize max}}(\hat{\mathbf{F}}_k+\hat{\mathbf{F}}^T_k)}{(a+h)\|\hat{\mathbf{P}}_k\|+c\|\hat{\mathbf{W}}_k+\hat{\mathbf{Y}}_k\|}\nonumber\\
=&~-\frac{1}{2}\frac{\lambda_{\mbox{\footnotesize max}}(\hat{\mathbf{F}}_k+\hat{\mathbf{F}}^T_k)}{(a+h)\|\hat{\mathbf{P}}_k\|+c\|\hat{\mathbf{P}}_k\mathbf{M}_k+\hat{\mathbf{P}}_k\mathbf{O}_k\|}\nonumber\\
\geq &~\frac{1}{2\|\hat{\mathbf{P}}_0\|}\frac{\gamma}{a+h+c\|\mathbf{M}_k+\mathbf{O}_k\|_F}\nonumber\\
\geq &~\frac{1}{2\|\hat{\mathbf{P}}_0\|}\frac{\gamma}{a+h+c\left(\mu+\omega\right)}.
\end{align}
\end{IEEEproof}

\section{Sparse Control Network Design}\label{sec:design}

To design a sparse control network, we seek a set of $\mathbf{L}_{k,ij}$ and $\mathbf{O}_{k,ij}$ that guarantee stability, with a small number links. Note that if link $ij$ is used it can carry both $\mathbf{L}_{k,ij}\hat{\mathbf{x}}_j(t)$ and $\mathbf{O}_{k,ij}\mathbf{y}_j(t)$. Thus, the use of a link from node $i$ to node $j$ at time $k$ can be encoded in binary variables $\alpha_{k,ij}\in\{0,1\}$. Then we can write $\mathbf{L}_{k,ij}=\alpha_{k,ij}\mathbf{L}^{\star}_{k,ij}$ and $\mathbf{O}_{k,ij}=\alpha_{k,ij}\mathbf{O}^{\star}_{k,ij}$, where $\mathbf{L}^{\star}_{k,ij}$ and $\mathbf{O}^{\star}_{k,ij}$ are the optimal link gains, when the link is used. In aggregate, these can be written as $\mathbf{L}_k=\boldsymbol{\alpha}_k\circ\mathbf{L}^{\star}_k$ and $\mathbf{O}_k=\hat{\boldsymbol\alpha}_k\circ\mathbf{O}^{\star}_k$, where $\boldsymbol\alpha_k=\left[\alpha_{k,ij}\mathbf{1}^{m_i\times n_j}\right]$ and $\hat{\boldsymbol\alpha}_k=\left[\alpha_{k,ij}\mathbf{1}^{n_i\times r_j}\right]$ with $\alpha_{k,ii}\equiv 0$.

Now, with the stability conditions provided in (\ref{eq: 39}) in hand, our objective is to design a control network with minimum number of links that satisfies stability conditions (\ref{eq: 39}). Minimizing the number of communication links is equivalent to minimizing the number of $\alpha_{k,ij}=1$, or in other words, minimizing the sum of $\alpha_{k,ij}$ subject to constraints in (\ref{eq: 39}). Our problem can therefore be formulated as the following convex mixed-binary program:
\begin{subequations}
\label{eq: 39.11}
\begin{align}
\text{minimize}&~\sum_{i,j\in \mathcal{N}}~\alpha_{k,ij}\label{eq: 39.11a}\\
\mbox{subject to}  &~\mbox{(\ref{eq: 39})}\label{eq: 39.11b}\\
&~\alpha_{k,ij}\in \{0,1\}\label{eq: 39.11c}
\end{align}
\end{subequations}
where $\mathbf{Y}_k=\boldsymbol{\alpha}_k \circ\mathbf{X}_k$ and $\hat{\mathbf{Y}}_k=\hat{\boldsymbol\alpha}_k \circ \hat{\mathbf{X}}_k$.


The complexity of solving problem \eqref{eq: 39.11} is important, since it has to be solved in each iteration. In general a mixed-binary program is NP-hard. In the worst case, one has to solve $\mathcal{O}(2^{N^2})$ convex problems, carrying an exhaustive search on the binary variables. While a variety of exact methods for convex mixed-binary programs are available 
\cite{grossmann2002}, their computational complexity is prohibitive for large networks, specially since the calculation is to be repeated periodically. Here, we propose a simple suboptimal relaxation-thresholding approach which should be carried out in each iteration for all $i,j\in \mathcal{N}$ where $\mathbf{H}_{ij}(t)$ is not identically zero. Otherwise, $\alpha_{k,ij}=0$ for all $k\in \mathbb{Z}_{+}$ as controller network is always a subset of plant network \cite{Razeghi2015}:
\begin{enumerate}
\item{Set $k\leftarrow 0$ and  $t_0\leftarrow 0$.}
\item{Find a feasible point for (\ref{eq: 39.11}) excluding (\ref{eq: 39e}) and (\ref{eq: 39f}) to find the initial conditions $\mathbf{Z}_0$ and $\hat{\mathbf{P}}_0$.}
\item{Initialize $\alpha_{k,ij}\leftarrow 1$ for all $i,j\in \mathcal{N},i\neq j$.}
\item{Find a feasible point for (\ref{eq: 39.11}) to yield $\mathbf{Z}_{k,i}$,$\mathbf{W}_{k,i}$,$\mathbf{X}_{k,ij}$,$\hat{\mathbf{P}}_{k,i}$,$\hat{\mathbf{W}}_{k,i}$ and $\hat{\mathbf{X}}_{k,ij}$. If a feasible solution is found, $\alpha^{\dagger}_{k,ij}\leftarrow \alpha_{k,ij}$. Otherwise go to step 7, unless the problem is infeasible at the first iteration, in which case there is no solution and the design procedure is terminated.}
\item{Solve \eqref{eq: 39.11} with \eqref{eq: 39.11c} relaxed to $\alpha_{k,ij}\in [0,1]$ to obtain solution $\alpha^{(r)}_{k,ij}$ satisfying \eqref{eq: 39a} and \eqref{eq: 39b} and $\mathbf{Z}_{k,i}$,$\mathbf{W}_{k,i}$,$\mathbf{X}_{k,ij}$,$\hat{\mathbf{P}}_{k,i}$,$\hat{\mathbf{W}}_{k,i}$,$\hat{\mathbf{X}}_{k,ij}$ are those found in step 4.}
\item{If all $\alpha^{(r)}_{k,ij}=0$, set $\alpha^{\dagger}_{k,ij}\leftarrow 0$ and go to step 7. Otherwise, set $\alpha_{k,ij}$ corresponding to the smallest non-zero $\alpha^{(r)}_{k,ij}$ to zero and return to step 4.}
\item{Return $\alpha^{\dagger}_{k,ij}$.}
\item{Calculate the next switching time $T_k$ from (\ref{eq: 39.14}).}
\end{enumerate}
Note that with the above design procedure, in the worst case, one has to solve $\mathcal{O}(N^2)$ convex problems in each iteration, since it can be solved by using a linear search in a sorted set $\{\alpha_{k,ij}^{(r)}\}$.

To further simplify the procedure we can substitute steps 6 and 7 above with
\begin{enumerate}
\setcounter{enumi}{5}
\item{Solve
\begin{eqnarray}
\label{eq:tau}
\begin{array}{rl}
\underset{m,l}{\text{maximize}}& \tau_k=\alpha^{(r)}_{k,ml}\\
\mbox{subject to}&  \alpha_{k,ij}=1_{\alpha^{(r)}_{k,ij}\geq \tau_k}, \mbox{\eqref{eq: 39a} and \eqref{eq: 39b}}.
\end{array}
\end{eqnarray}}
\item{Return $\alpha^{\dagger}_{k,ij}=1_{\alpha_{k,ij}^{(r)}>\tau^{\star}_k}$, where $\tau^\star_k$ is the solution of \eqref{eq:tau}.}
\end{enumerate}
We note that the maximum number of convex problems that should be solved in \eqref{eq:tau} in each iteration is only $\mathcal{O}(\log N)$, since it can be solved by a binary search on $\tau_k$ in a sorted set $\{\alpha_{k,ij}^{(r)}\}$. Of course, this reduction in complexity is at the price of a more conservative solution.\

In the next section, we will show the gap between this suboptimal relaxation-thresholding approach and the optimal exhaustive search and how our heuristic algorithm performs near optimal.

\section{Numerical Example}\label{sec:example}

\begin{figure}[!t]
\centering
\includegraphics[width =2.8in]{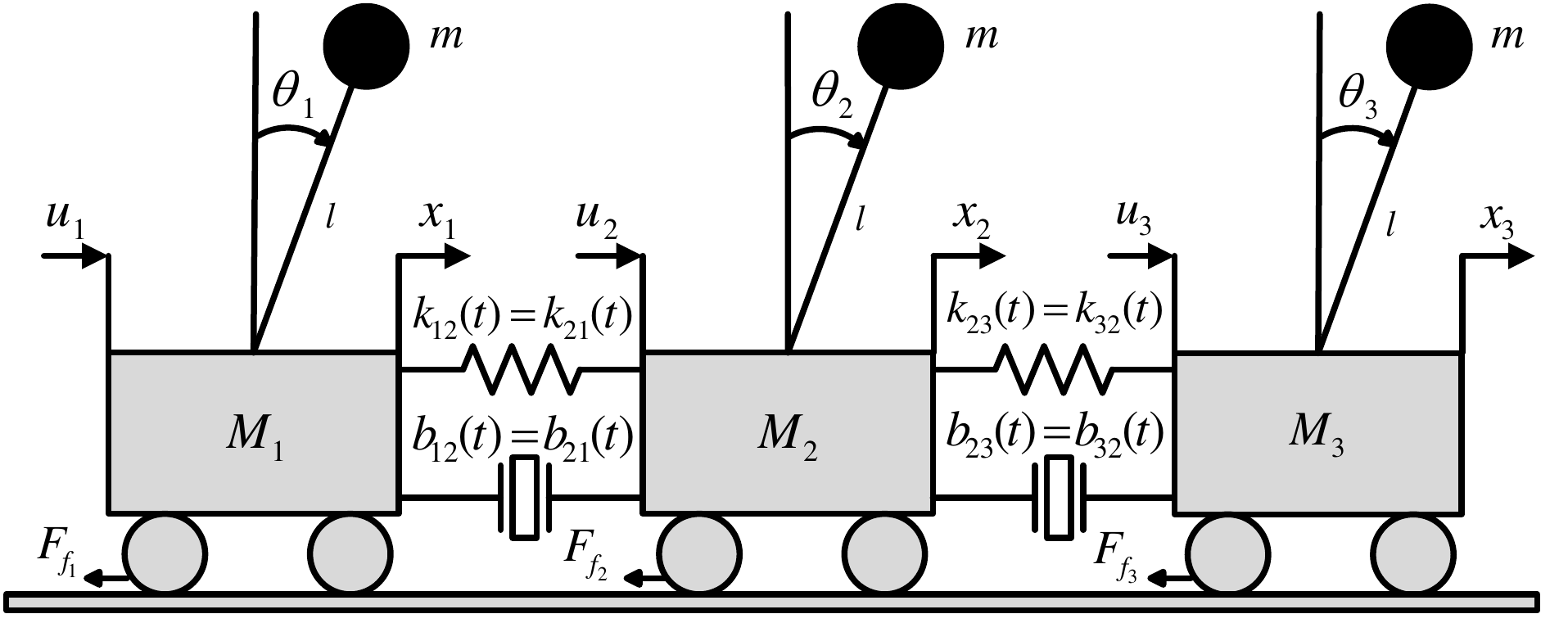}
\vspace{-0.1in}
\caption {Network of three coupled inverted pendulums}
\vspace{-0.2in}
\label{fig:NCS2}
\end{figure}

Consider the system shown in Fig. \ref{fig:NCS2}, where three inverted pendulums are mounted on coupled carts. Linearizing equations of motions yield \cite{ogata}
\begin{eqnarray}
\label{eq: 302}
\mathbf{A}_i(t)&=&\left[
               \begin{array}{cccc}
                 0 & 1 & 0 & 0 \\
                 \frac{M_i+m}{M_il}g & 0 & \frac{k_i(t)}{M_il} & \frac{c_i+b_i(t)}{M_il} \\
                 0 & 0 & 0 & 1 \\
                 \frac{-m}{M_i}g & 0 & \frac{-k_i(t)}{M_i} & \frac{-c_i-b_i(t)}{M_i}
               \end{array}
             \right],\nonumber\\
\mathbf{B}_i(t)&=&\left[\begin{array}{cccc} 0 & \frac{-1}{M_il} & 0 & \frac{1}{M_i}\end{array}\right]^T,~\mathbf{C}_i(t)=\left[\begin{array}{cccc}1 & 0 & 0 & 0 \\0 & 0 & 1 & 0\end{array}\right],\nonumber\\
\mathbf{H}_{ij}(t)&=&\left[
          \begin{array}{cccc}
            0 & 0 & 0 & 0 \\
            0 & 0 & \frac{-k_{ij}(t)}{M_il} & \frac{-b_{ij}(t)}{M_il} \\
            0 & 0 & 0 & 0 \\
            0 & 0 & \frac{k_{ij}(t)}{M_i} & \frac{b_{ij}(t)}{M_i}
          \end{array}
        \right],
\end{eqnarray}
for $(i,j)=(1,2),(2,1),(2,3),(3,2)$, where $k_i(t)=\sum_{j\in \mathcal{N}_i}k_{ij}(t)$ and $b_i(t)=\sum_{j\in \mathcal{N}_i}b_{ij}(t)$. Here $c_i$, $b_{ij}(t)=b_{ji}(t)$, and $k_{ij}(t)=k_{ji}(t)$ are friction, damper and spring coefficients, respectively, and we have assumed the moment of inertia of the pendulums to be zero.\

Since the subsystems are controllable and observable, we can use the optimization problem \eqref{eq: 39.11} to design distributed observers and controllers that stabilize the entire network with small number of links in the control network. As design criteria, we assume bounds on the norm of local gains are $\kappa_1=\kappa_2=280,~\kappa_3=480$ and $\mu_i=40$ and bounds on the norms of coupling gains are $\iota_{ij}=20$ and $\omega_{ij}=10$ and the numerical system parameters are $M_1=5$, $M_2=3$, $M_3=7$, $m=1$, $g=10$, $l=1$, $k_{ij}(t)=1+0.5\cos(t)$, $b_{ij}(t)=1+0.5\sin(t)$, $c_1=4$, $c_2=2$ and $c_3=1$.

All subsystem matrices are continuously differentiable uniformly bounded with $a=0.48$, $h=0.34$ and $b=c=0$. Hence, Assumption 1 is satisfied. We set $\gamma=0.2$, $\beta_i=0.01$ and $\epsilon_i=0.05$ for all $i$.\

\begin{figure}[!t]
\centering
\includegraphics[width =2.8in]{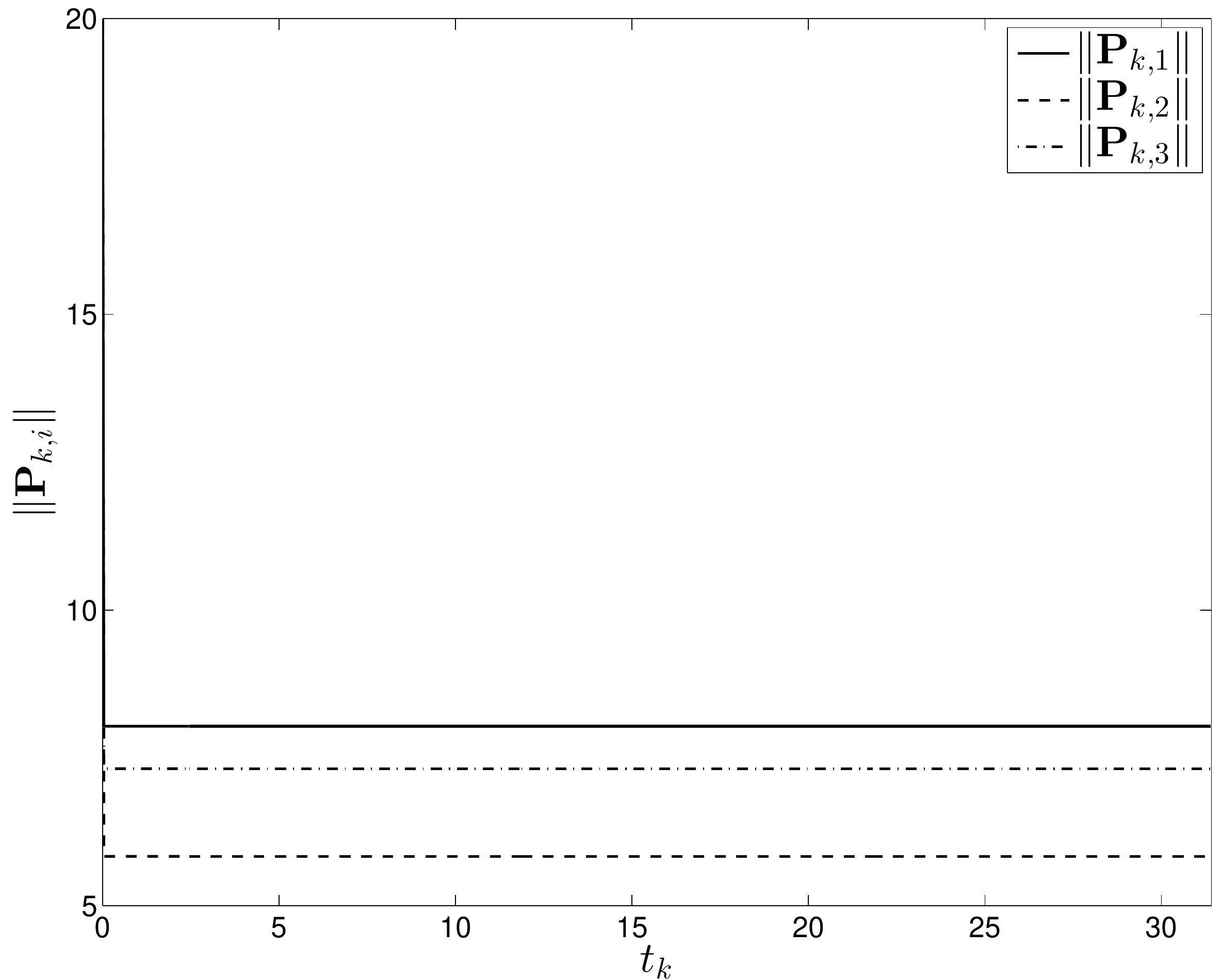}
\caption {Norm of $\mathbf{P}_{k,i}$ versus time.}
\label{fig:P}
\end{figure}

\begin{figure}[!t]
\centering
\includegraphics[width =2.8in]{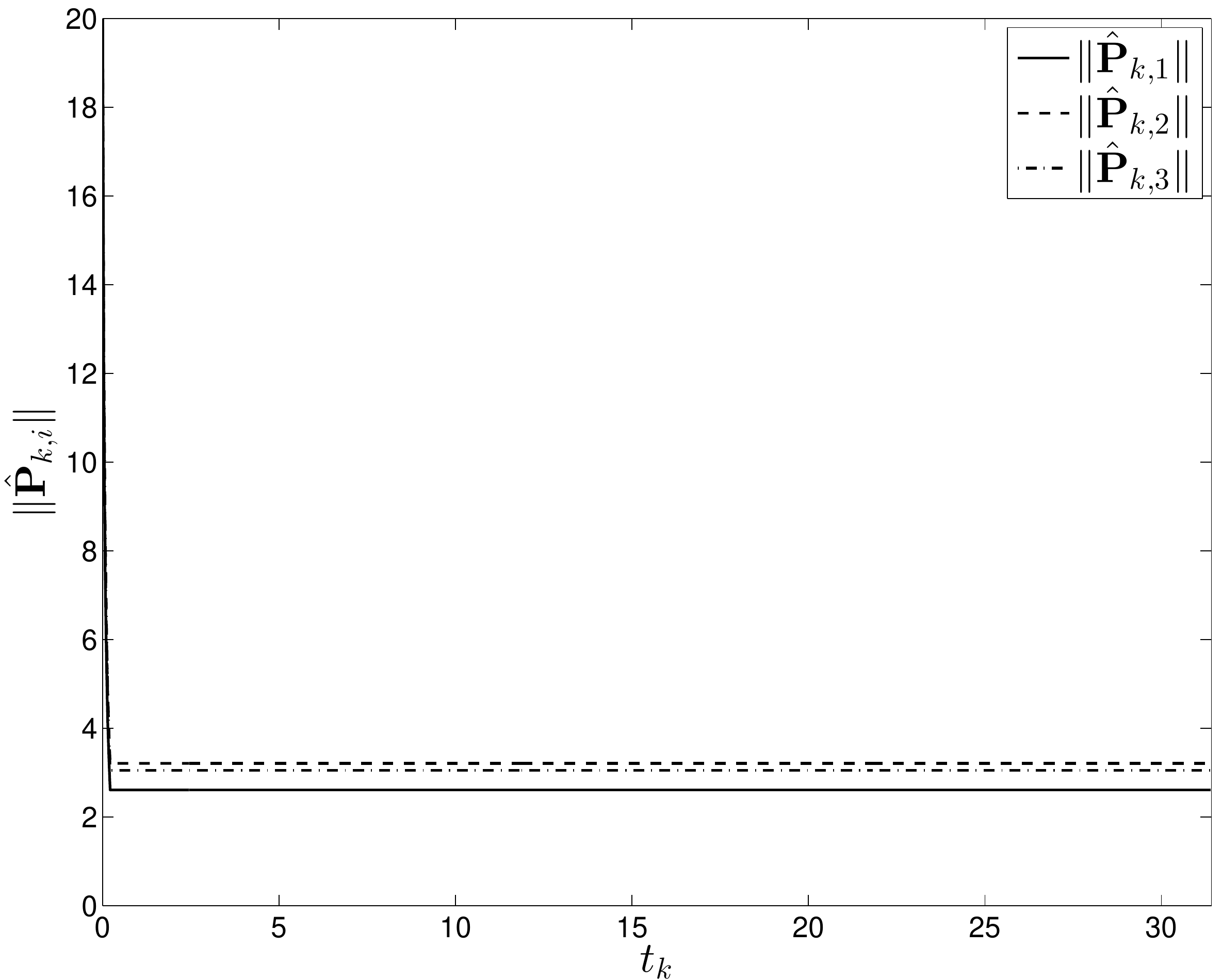}
\caption {Norm of $\hat{\mathbf{P}}_{k,i}$ versus time.}
\label{fig:Phat}
\end{figure}

\begin{figure}[!t]
\centering
\includegraphics[width =2.8in]{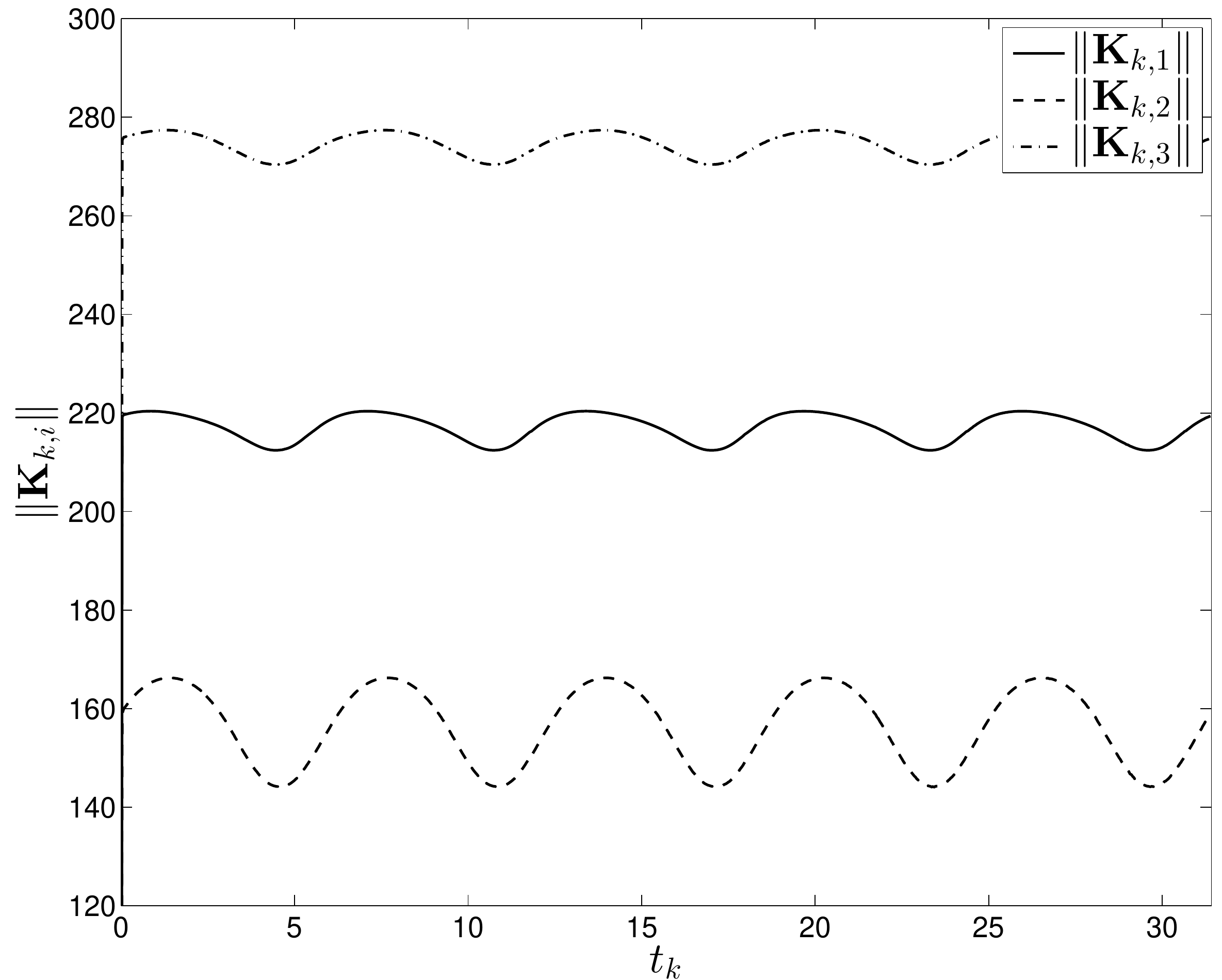}
\caption {Norm of local controller gains $\mathbf{K}_{k,i}$ versus time.}
\label{fig:K}
\end{figure}

\begin{figure}[!t]
\centering
\includegraphics[width =2.8in]{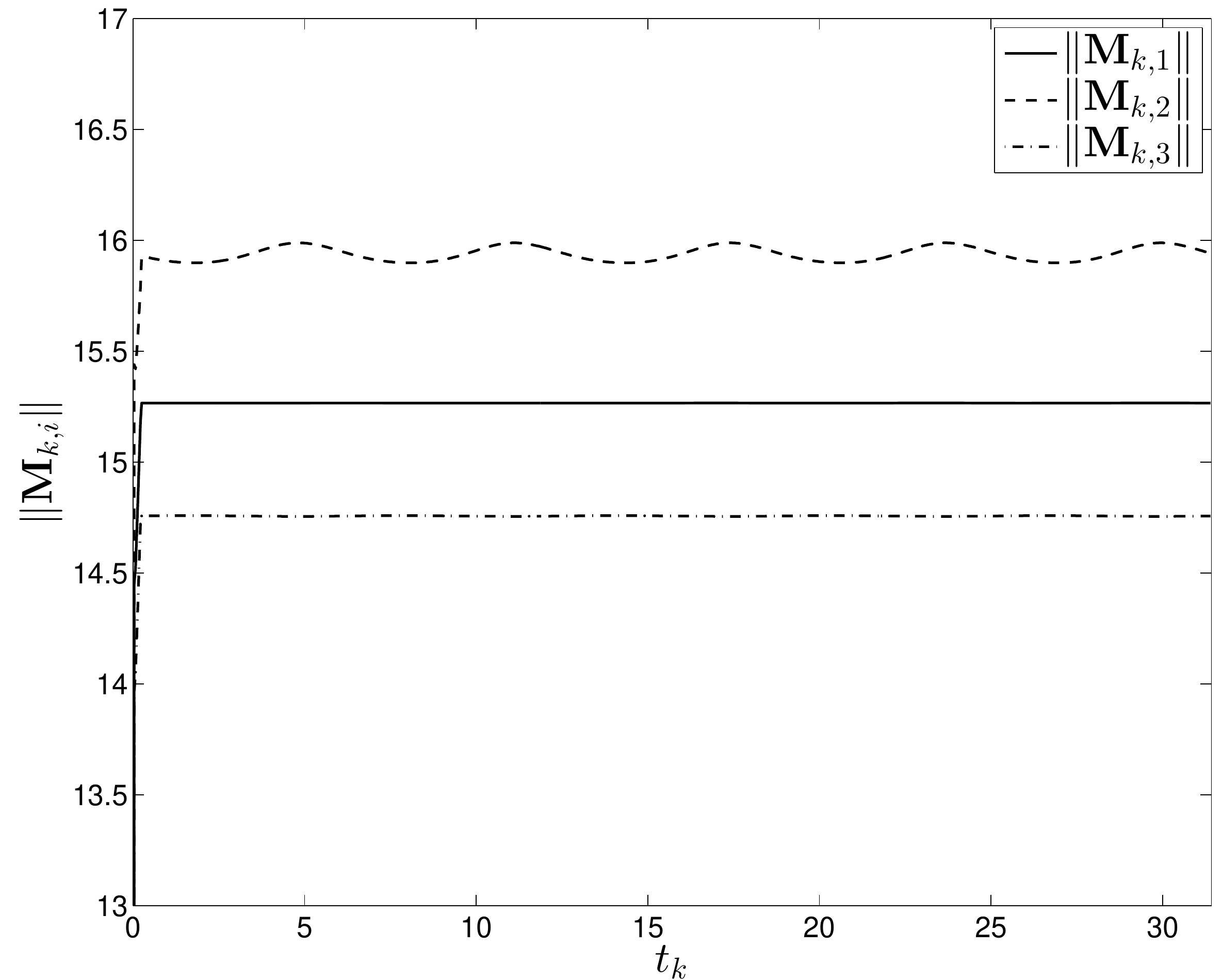}
\caption {Norm of local controller gains $\mathbf{M}_{k,i}$ versus time.}
\label{fig:M}
\end{figure}

\begin{figure}[!t]
\centering
\includegraphics[width =2.8in]{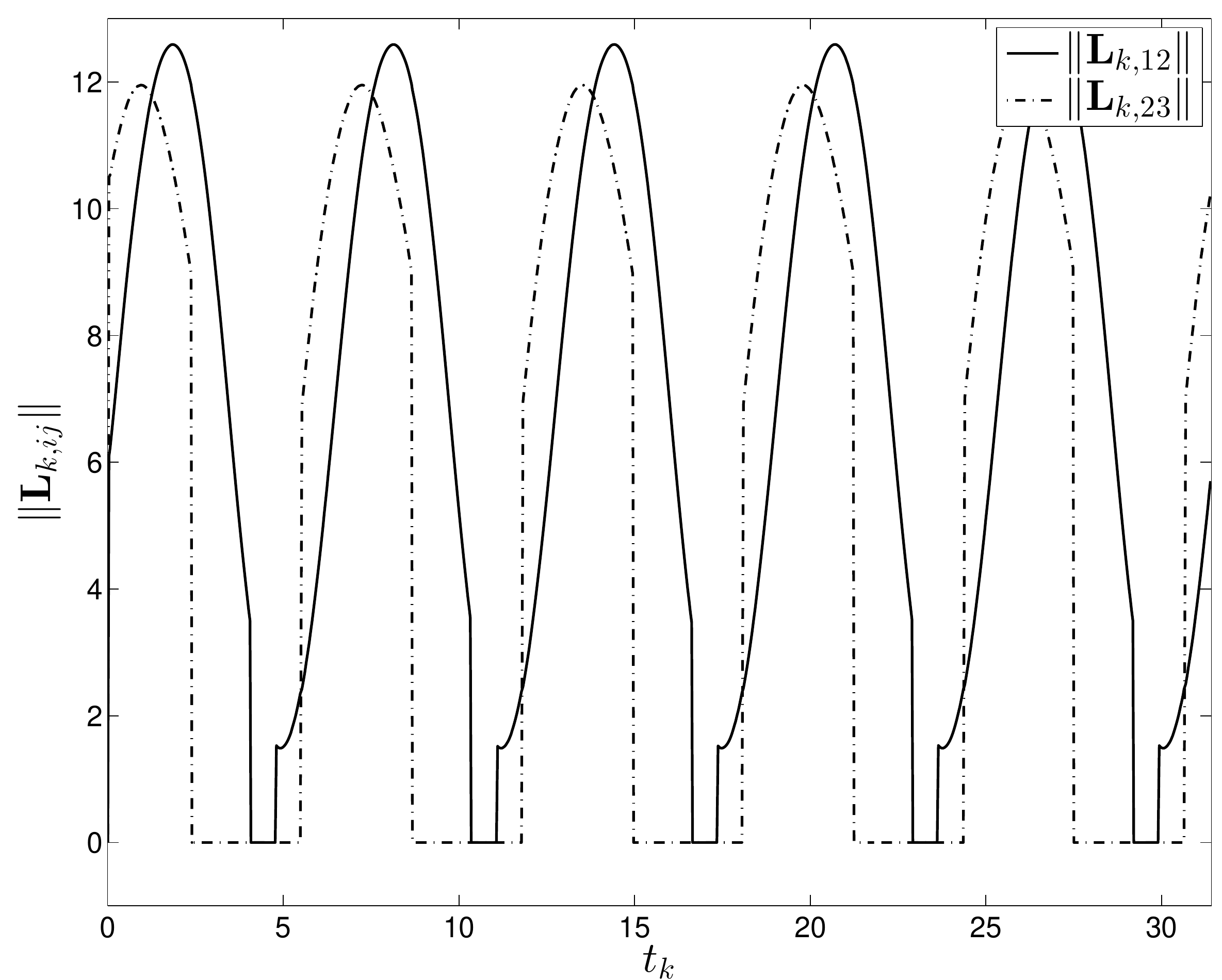}
\caption {Norm of coupling controller gains $\mathbf{L}_{k,ij}$ versus time. Note that the rest of links, i.e., $\mathbf{L}_{k,21}=\mathbf{L}_{k,32}=\mathbf{L}_{k,13}=\mathbf{L}_{k,31}=0$ for all $t_k$.}
\label{fig:L}
\end{figure}

\begin{figure}[!t]
\centering
\includegraphics[width =2.8in]{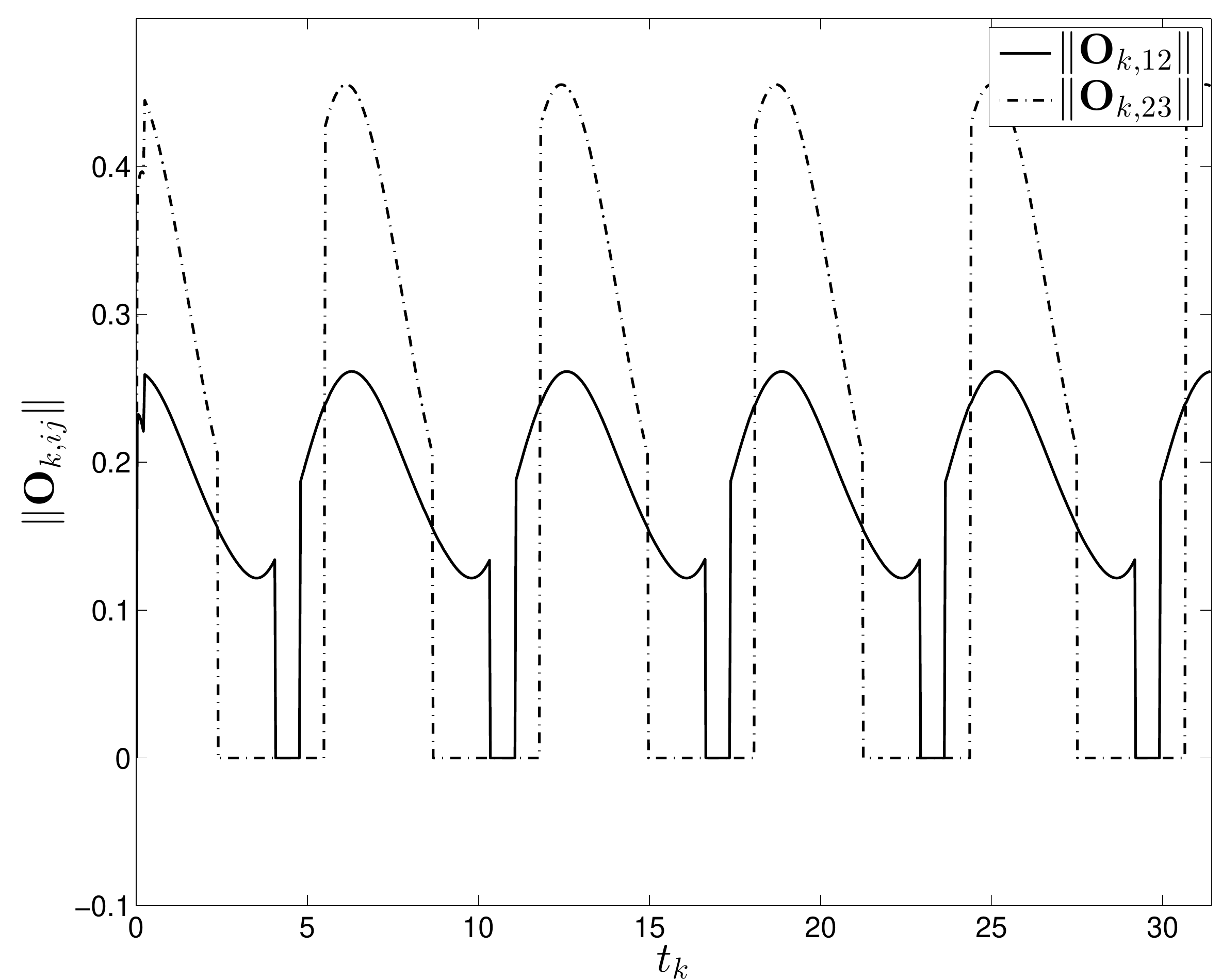}
\caption {Norm of coupling observer gains $\mathbf{O}_{k,ij}$ versus time. Note that the rest of links, i.e., $\mathbf{O}_{k,21}=\mathbf{O}_{k,32}=\mathbf{O}_{k,13}=\mathbf{O}_{k,31}=0$ for all $t_k$.}
\label{fig:O}
\end{figure}

\begin{figure}[!t]
\centering
\includegraphics[width =2.8in]{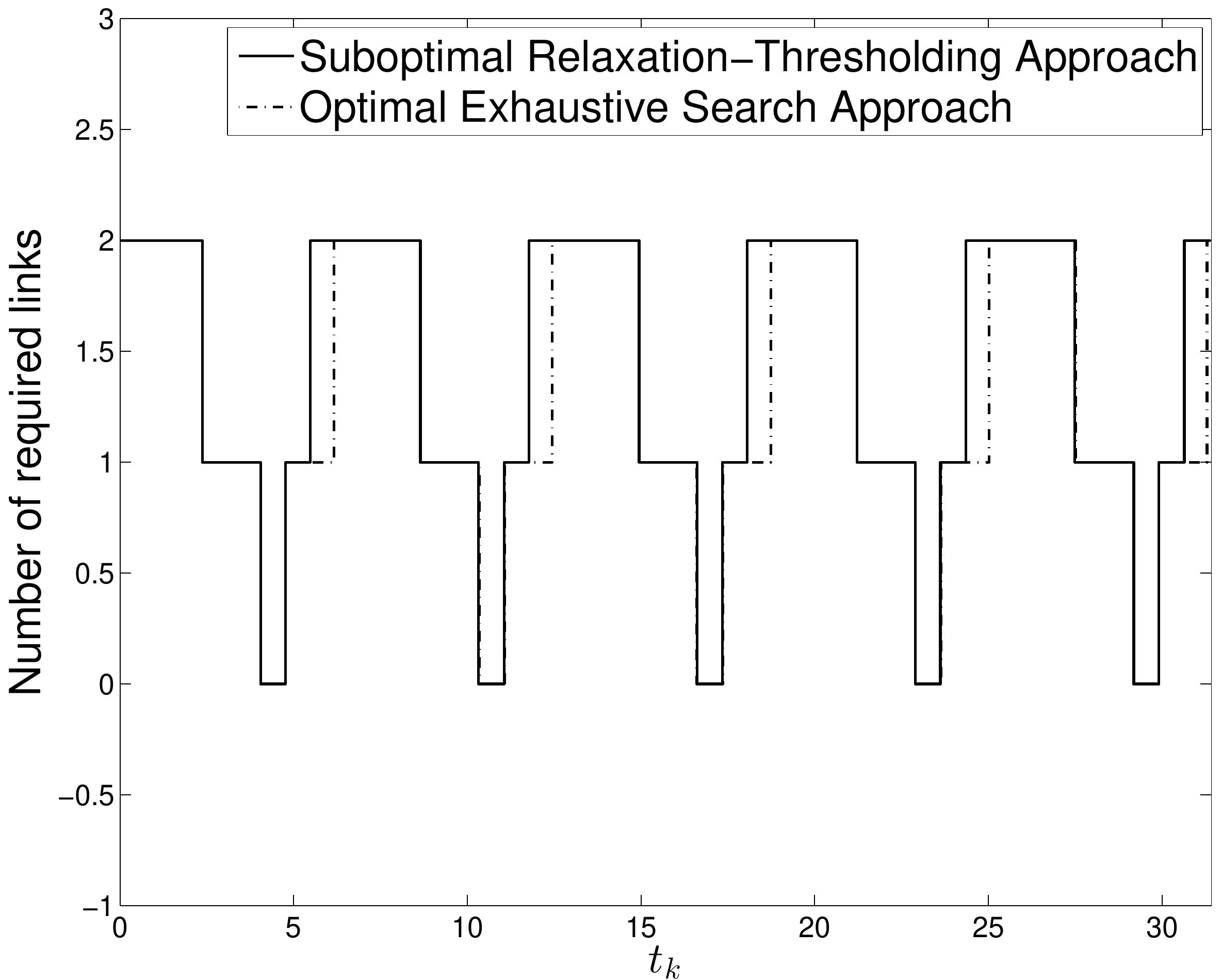}
\caption {Number of required links versus time for the two cases: (i) proposed simple suboptimal relaxation-thresholding approach and (ii) optimal exhaustive search approach.}
\label{fig:Links}
\end{figure}

\begin{figure}[!t]
\centering
\includegraphics[width =2.8in]{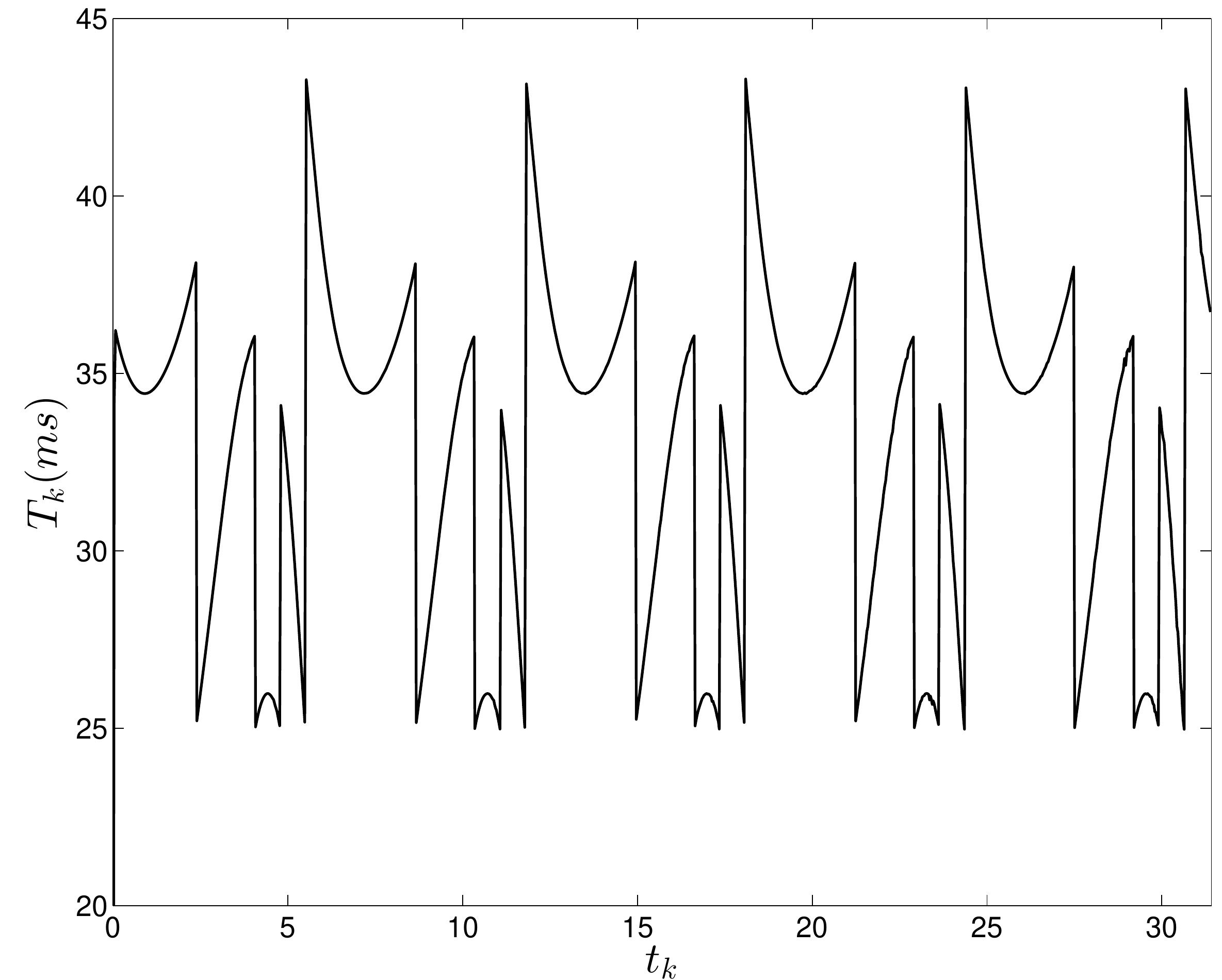}
\caption {Updating times, $T_k$.}
\label{fig:T}
\end{figure}

The simulation results are presented in Fig. \ref{fig:P} to Fig. \ref{fig:T} for $t=[0,10\pi]$ seconds where to solve (\ref{eq: 39.11}), we used our proposed simple suboptimal relaxation-thresholding approach with linear search. To show that how this performs well (near optimal), we compare the number of required links with the optimal exhaustive search on binary variables. Fig. \ref{fig:P} and \ref{fig:Phat} depict $\|\mathbf{P}_{k,i}\|$ and $\|\hat{\mathbf{P}}_{k,i}\|$ as a function of $t_k$. We can see that, as expected they converge as $t_k\rightarrow \infty$. The convergence of $\|\mathbf{P}_{k,i}\|$ happens after only one time slot, while $\|\hat{\mathbf{P}}_{k,i}\|$ takes 5 time slots to converge. Fig. \ref{fig:K} and Fig. \ref{fig:M}, depict $\|\mathbf{K}_{k,i}\|$ and $\|\mathbf{M}_{k,i}\|$, which are the local controller and observer gains, respectively. Similarly, Fig. \ref{fig:L}, and Fig. \ref{fig:O} depict $\|\mathbf{L}_{k,ij}\|$ and $\|\mathbf{O}_{k,ij}\|$, which are the coupling controller and observer gains, respectively. We observe that whenever a link is not necessary, (i.e., $\alpha_{k,ij}=0$), the link gain is set to zero. Otherwise, it is assigned the optimal values $\mathbf{L}^{\star}_{k,ij}$ and $\mathbf{O}^{\star}_{k,ij}$. We can see that all local and coupling gains are limited as enforced in \eqref{eq: 5}.

In Fig. \ref{fig:Links}, we plot number of required links in communications network versus $t_k$ for the two cases: (i) proposed simple suboptimal relaxation-thresholding approach and (ii) optimal exhaustive search on binary variables. The suboptimal approach performs very well and as we expected, the optimal search gives better result at the price of more complexity. Fig. \ref{fig:T} presents the updating times, $T_k$ versus $t_k$. We see that the fewer number of links indicates the shorter updating time $T_k$ and whenever there is a change in the number of required links, there also is a step change in the updating time $T_k$.

\section{Concluding Remarks}
\label{sec:conclusion}
We have provided an iterative design approach for distributed observer-based controllers that stabilize a given linear time-varying networked control system with an arbitrary directed topology. To measure states of each subsystem, we use the outputs of other subsystems to improve stability of observer dynamics; this approach is the dual of the distributed controller network. Our design approach is based on a set of stability conditions obtained using the piecewise quadratic switching stabilization method with multiple quadratic Lyapunov functions, which must be updated and optimized in discrete time and provides a sparse observer-controller network that guarantees global asymptotic stability. Due to the assumptions made here to maintain tractability, the design has some degree of conservatism. Thus, although the results provide us with significant insight into the problem of designing a sparse observer-controller network, a gap still remains. Further quantification or reduction of this gap will be quite valuable.

We added a free variable to the stability inequalities to avoid spending the entire margin in the stability criteria during the search for a sparse observer-controller network. Optimal distribution of this margin among the inequalities to make the network robust without significantly growing the size of the observer-controller network is, however, unknown. Therefore, further investigation of the tradeoff between the stability margin and the sparsity of the observer-controller network will be interesting.

We believe that the results presented in this paper provide a foundation for further progress toward understanding these interesting and important problems.

\section*{Acknowledgment}
The authors would like to give special thanks to Prof. Mark F. Bocko at University of Rochester for his support and helpful suggestions.


\vspace{-.5in}
\begin{IEEEbiography}[{\includegraphics[width=1in,height=1.25in,clip,keepaspectratio]{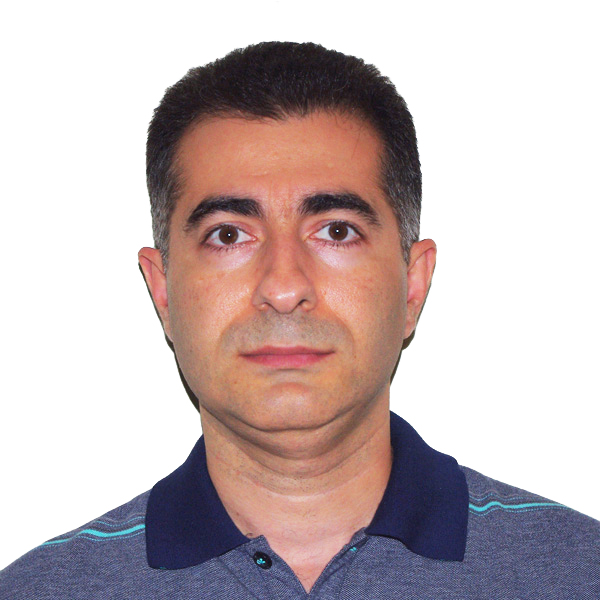}}]{Mohammad Razeghi-Jahromi} received his B.S. degree from Amirkabir University of Technology, Tehran, Iran, in 1997 and his M.S. degrees from University of Tehran, Tehran, Iran in 2000, and his Ph.D. degree from University of Rochester, Rochester, NY, USA in 2016, all in Electrical Engineering. He is now research scientist at ABB corporate research United States (USCRC), Raleigh, NC, USA, since 2017. His primary research interests include networked control systems, control systems theory, stochastic control and stochastic differential equations, Markov jump linear systems, and convex optimization.
\end{IEEEbiography}

\vspace{-.5in}
\begin{IEEEbiography}[{\includegraphics[width=1in,height=1.25in,clip,keepaspectratio]{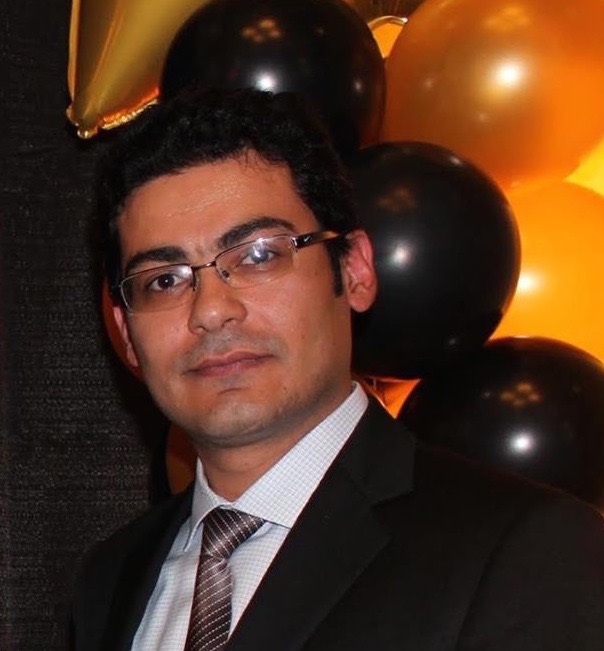}}]{Saeed Manaffam}
received his B.S. degree in Electrical Engineering from Amirkabir University of Technology, Tehran, Iran, in 2006 and his M.Sc. degrees also in Electrical Engineering from Sharif University of Technology, Tehran, Iran in 2008 and University of Rochester, NY, USA in 2012, respectively. He is now working toward his Ph.D. degree in the department of Electrical Engineering and Computer Science at University of Central Florida, Orlando, FL, USA, since August 2012. His current research includes the study of dynamical systems in complex networks, decentralized and distributed control of such networks.
\end{IEEEbiography}

\vspace{-.5in}
\begin{IEEEbiography}[{\includegraphics[width=1in,height=1.25in,clip,keepaspectratio]{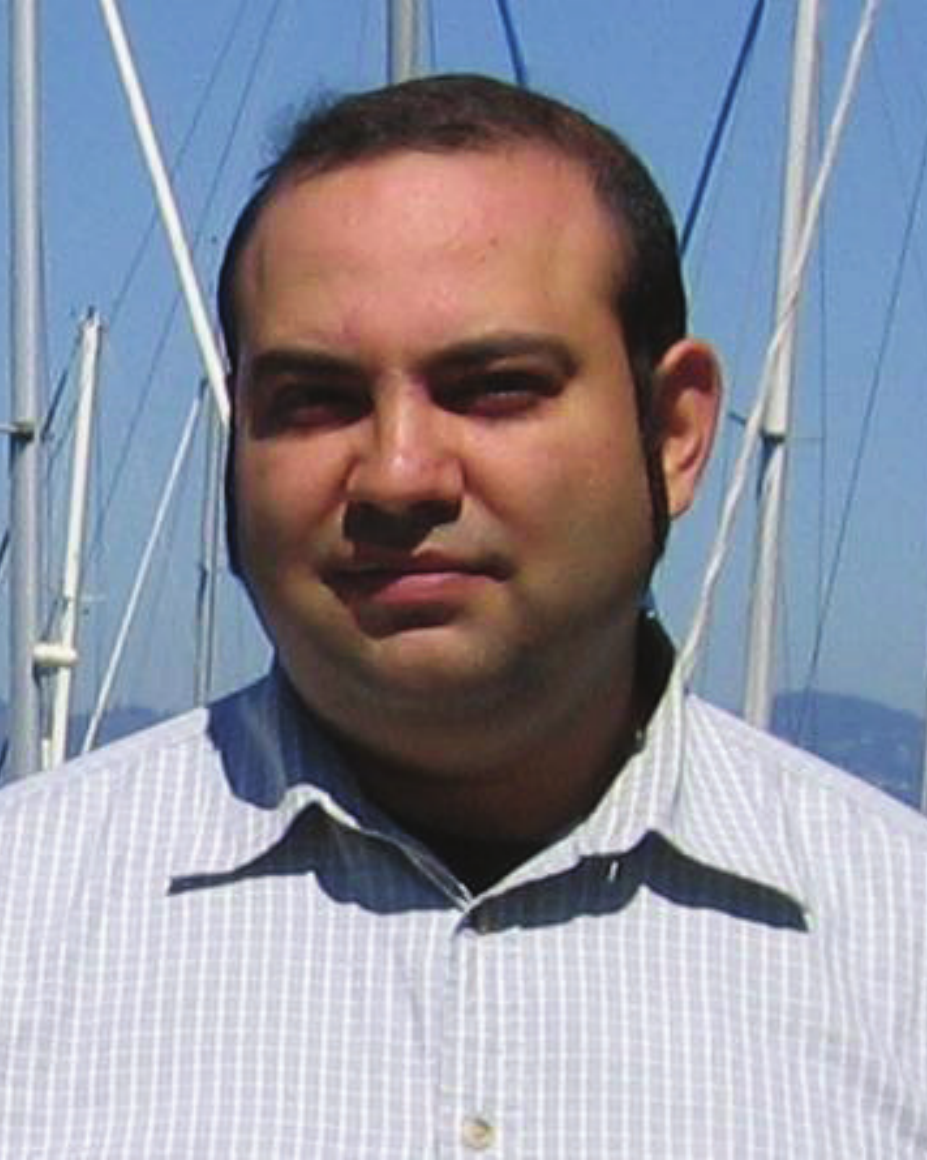}}]{Alireza Seyedi}
(S'95-M'04-SM'10) was a Research Assistant Professor in the Department of Electrical Engineering and Computer Science at the University of Central Florida. He had been with the department since August of 2012. Prior to that, he had been at the University of Rochester and at Philips Research North America. He had received his Ph.D. and M.S. degrees from Rensselaer Polytechnic Institute, both in Electrical Engineering, in 2004 and 2000, respectively. He had received his B.S. degree, also in Electrical Engineering, from Sharif University of Technology in 1997. His research interests were in the convergence of communications and control. In particular, he worked on control and decision making for networked systems, distributed and decentralized control, dynamics and stochastics of complex networks, control and communication for energy networks and the Smart Grid, decision and control for systems with stochastic sources of energy, energy harvesting for communications and cognitive radios and networks. He was a member of Eta Kappa Nu and Tau Beta Pi and was an Associate Editor of the \sc{IEEE Signal Processing Letters}.
\end{IEEEbiography}

\vspace{-.5in}
\begin{IEEEbiography}[{\includegraphics[width=1in,height=1.25in,clip,keepaspectratio]{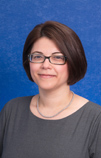}}]{Azadeh Vosoughi}
(S'95-M'06-SM'14) Azadeh Vosoughi is an Associate Professor in the Department of Electrical Engineering and Computer Science at UCF. She received her B.S. degree from Sharif University of Technology, Tehran, Iran, in 1997, her M.S. degree from Worcester Polytechnic Institute, Worcester, MA, in 2001, and her Ph.D. degree from Cornell University, Ithaca, NY, in 2006, all in Electrical Engineering. She received the NSF CAREER award in 2011 for her research on the integration of signal processing and communications for decision fusion in distributed detection systems. Prior to joining the faculty at UCF, she was an Assistant Professor at the University of Rochester. Currently, she is serving as an Associate Editor for the {\sc{IEEE Transactions on Wireless Communications}} and the {\sc{IEEE Signal Processing Letters}}.
\end{IEEEbiography}


\begin{thebibliography}{10}
	\providecommand{\url}[1]{#1}
	\csname url@samestyle\endcsname
	\providecommand{\newblock}{\relax}
	\providecommand{\bibinfo}[2]{#2}
	\providecommand{\BIBentrySTDinterwordspacing}{\spaceskip=0pt\relax}
	\providecommand{\BIBentryALTinterwordstretchfactor}{4}
	\providecommand{\BIBentryALTinterwordspacing}{\spaceskip=\fontdimen2\font plus
		\BIBentryALTinterwordstretchfactor\fontdimen3\font minus
		\fontdimen4\font\relax}
	\providecommand{\BIBforeignlanguage}[2]{{%
			\expandafter\ifx\csname l@#1\endcsname\relax
			\typeout{** WARNING: IEEEtran.bst: No hyphenation pattern has been}%
			\typeout{** loaded for the language `#1'. Using the pattern for}%
			\typeout{** the default language instead.}%
			\else
			\language=\csname l@#1\endcsname
			\fi
			#2}}
	\providecommand{\BIBdecl}{\relax}
	\BIBdecl
	
	\bibitem{Siljak1991}
	D.~Siljak, \emph{Decentralized control of complex systems}.\hskip 1em plus
	0.5em minus 0.4em\relax Dover Publications, 2012.
	
	\bibitem{Siljak2005}
	A.~Zecevic and D.~Siljak, ``Global low-rank enhancement of decentralized
	control for large-scale systems,'' \emph{IEEE Transactions on Automatic
		Control}, vol.~50, no.~5, pp. 740--744, 2005.
	
	\bibitem{Mazo2010}
	M.~Mazo, A.~Anta, and P.~Tabuada, ``An {ISS} self-triggered implementation of
	linear controllers,'' \emph{Automatica}, vol.~46, no.~8, pp. 1310--1314,
	2010.
	
	\bibitem{Zhang2001}
	W.~Zhang, M.~Branicky, and S.~Phillips, ``Stability of networked control
	systems,'' \emph{IEEE Control Systems}, vol.~21, no.~1, pp. 84--99, 2001.
	
	\bibitem{Baillieul2007}
	J.~Baillieul and P.~Antsaklis, ``Control and communication challenges in
	networked real-time systems,'' \emph{Proceedings of the IEEE}, vol.~95,
	no.~1, pp. 9--28, 2007.
	
	\bibitem{Hespanha2007}
	J.~Hespanha, P.~Naghshtabrizi, and Y.~Xu, ``A survey of recent results in
	networked control systems,'' \emph{Proceedings of the IEEE}, vol.~95, no.~1,
	pp. 138--162, 2007.
	
	\bibitem{Yue2008}
	F.~Wang and D.~Liu, \emph{Networked control systems: theory and
		applications}.\hskip 1em plus 0.5em minus 0.4em\relax Springer, 2008.
	
	\bibitem{Wang2011}
	X.~Wang and M.~Lemmon, ``Event-triggering in distributed networked control
	systems,'' \emph{IEEE Transactions on Automatic Control}, vol.~56, no.~3, pp.
	586--601, 2011.
	
	\bibitem{montestruque2004}
	L.~Montestruque and P.~Antsaklis, ``Stability of model-based networked control
	systems with time-varying transmission times,'' \emph{IEEE Transactions on
		Automatic Control}, vol.~49, no.~9, pp. 1562--1572, 2004.
	
	\bibitem{Teel2004}
	D.~Nesic and A.~Teel, ``Input-output stability properties of networked control
	systems,'' \emph{IEEE Transactions on Automatic Control}, vol.~49, no.~10,
	pp. 1650--1667, 2004.
	
	\bibitem{Walsh2002}
	G.~Walsh, H.~Ye, and L.~Bushnell, ``Stability analysis of networked control
	systems,'' \emph{IEEE Transactions on Control Systems Technology}, vol.~10,
	no.~3, pp. 438--446, 2002.
	
	\bibitem{seiler2005h}
	P.~Seiler and R.~Sengupta, ``An $\mathcal{H}_{\infty}$ approach to networked
	control,'' \emph{IEEE Transactions on Automatic Control}, vol.~50, no.~3, pp.
	356--364, 2005.
	
	\bibitem{yu2005}
	M.~Yu, L.~Wang, T.~Chu, and F.~Hao, ``Stabilization of networked control
	systems with data packet dropout and transmission delays: continuous-time
	case,'' \emph{European Journal of Control}, vol.~11, no.~1, pp. 40--49, 2005.
	
	\bibitem{wu2007}
	J.~Wu and T.~Chen, ``Design of networked control systems with packet
	dropouts,'' \emph{IEEE Transactions on Automatic Control}, vol.~52, no.~7,
	pp. 1314--1319, 2007.
	
	\bibitem{witrant2007}
	E.~Witrant, C.~Canudas-de Wit, D.~Georges, and M.~Alamir, ``Remote
	stabilization via communication networks with a distributed control law,''
	\emph{IEEE Transactions on Automatic Control}, vol.~52, no.~8, pp.
	1480--1485, 2007.
	
	\bibitem{Gao08}
	H.~Gao, T.~Chen, and J.~Lam, ``A new delay system approach to network-based
	control,'' \emph{Automatica}, vol.~44, no.~1, pp. 39 -- 52, 2008.
	
	\bibitem{Zhang13}
	L.~Zhang, H.~Gao, and O.~Kaynak, ``Network-induced constraints in networked
	control systems: A survey,'' \emph{IEEE Trans. on Industrial Informatics},
	vol.~9, no.~1, pp. 403--416, Feb 2013.
	
	\bibitem{Lall2006}
	M.~Rotkowitz and S.~Lall, ``A characterization of convex problems in
	decentralized control,'' \emph{IEEE Transactions on Automatic Control},
	vol.~51, no.~2, pp. 274--286, 2006.
	
	\bibitem{Rotkowitz2009}
	M.~Rotkowitz and N.~Martins, ``On the closest quadratically invariant
	constraint,'' in \emph{Proceedings of the 48th IEEE Conference on Decision
		and Control}.\hskip 1em plus 0.5em minus 0.4em\relax IEEE, 2009, pp.
	1607--1612.
	
	\bibitem{Rotkowitz2010}
	M.~Rotkowitz, R.~Cogill, and S.~Lall, ``Convexity of optimal control over
	networks with delays and arbitrary topology,'' \emph{International Journal of
		Systems, Control and Communications}, vol.~2, no.~1, pp. 30--54, 2010.
	
	\bibitem{Swigart2009}
	J.~Swigart and S.~Lall, ``A graph-theoretic approach to distributed control
	over networks,'' in \emph{Proceedings of the 48th IEEE Conference on Decision
		and Control}.\hskip 1em plus 0.5em minus 0.4em\relax IEEE, 2009, pp.
	5409--5414.
	
	\bibitem{Shah2011}
	P.~Parrilo, P.~Shah \emph{et~al.}, ``A partial order approach to decentralized
	control,'' Ph.D. dissertation, Massachusetts Institute of Technology, 2011.
	
	\bibitem{Parrilo2011}
	P.~Shah and P.~Parrilo, ``An optimal controller architecture for poset-causal
	systems,'' in \emph{50th IEEE Conference on Decision and Control}.\hskip 1em
	plus 0.5em minus 0.4em\relax IEEE, 2011, pp. 5522--5528.
	
	\bibitem{Massoud2011}
	S.~Amin, ``Smart grid: Overview, issues and opportunities. advances and
	challenges in sensing, modeling, simulation, optimization and control,''
	\emph{European Journal of Control}, vol.~17, no. 5-6, pp. 547--567, 2011.
	
	\bibitem{Schuppen2011}
	J.~van Schuppen, O.~Boutin, P.~Kempker, J.~Komenda, T.~Masopust, N.~Pambakian,
	and A.~Ran, ``Control of distributed systems: Tutorial and overview,''
	\emph{European Journal of Control}, vol.~17, no.~5, pp. 579--602, 2011.
	
	\bibitem{Lafortune2003}
	K.~Rudie, S.~Lafortune, and F.~Lin, ``Minimal communication in a distributed
	discrete-event system,'' \emph{IEEE Transactions on Automatic Control},
	vol.~48, no.~6, pp. 957--975, 2003.
	
	\bibitem{Lin2007}
	F.~Lin, K.~Rudie, and S.~Lafortune, ``Minimal communication for essential
	transitions in a distributed discrete-event system,'' \emph{IEEE Transactions
		on Automatic Control}, vol.~52, no.~8, pp. 1495--1502, 2007.
	
	\bibitem{Razeghi2015}
	M.~Razeghi-Jahromi and A.~Seyedi, ``Stabilization of networked control systems
	with sparse observer-controller networks,'' \emph{IEEE Transactions on
		Automatic Control}, vol.~60, no.~6, pp. 1686--1691, 2015.
	
	\bibitem{liberzon1999basic}
	D.~Liberzon and A.~S. Morse, ``Basic problems in stability and design of
	switched systems,'' \emph{IEEE Control Systems}, vol.~19, no.~5, pp. 59--70,
	1999.
	
	\bibitem{liberzon2003switching}
	D.~Liberzon, \emph{Switching in systems and control}.\hskip 1em plus 0.5em
	minus 0.4em\relax Springer, 2003.
	
	\bibitem{michel1999recent}
	A.~N. Michel, ``Recent trends in the stability analysis of hybrid dynamical
	systems,'' \emph{IEEE Transactions on Circuits and Systems I: Fundamental
		Theory and Applications}, vol.~46, no.~1, pp. 120--134, 1999.
	
	\bibitem{lin2009stability}
	H.~Lin and P.~J. Antsaklis, ``Stability and stabilizability of switched linear
	systems: a survey of recent results,'' \emph{IEEE Transactions on Automatic
		control}, vol.~54, no.~2, pp. 308--322, 2009.
	
	\bibitem{haddad2011}
	W.~Haddad and V.~Chellaboina, \emph{Nonlinear dynamical systems and control: a
		Lyapunov-based approach}.\hskip 1em plus 0.5em minus 0.4em\relax Princeton
	University Press, 2011.
	
	\bibitem{grossmann2002}
	I.~Grossmann, ``Review of nonlinear mixed-integer and disjunctive programming
	techniques,'' \emph{Optimization and Engineering}, vol.~3, no.~3, pp.
	227--252, 2002.
	
	\bibitem{ogata}
	K.~Ogata, \emph{Modern control engineering, 1997}.\hskip 1em plus 0.5em minus
	0.4em\relax Prentice-Hall Inc., NJ.
	
\end{thebibliography}
\end{document}